\documentclass[11pt]{amsart}

 \newtheorem{thm}{Theorem}[section]
 \newtheorem{cor}[thm]{Corollary}
 \newtheorem{lem}[thm]{Lemma}
 \newtheorem{prop}[thm]{Proposition}
 \theoremstyle{definition}
 \newtheorem{defn}[thm]{Definition}
 \theoremstyle{remark}
 \newtheorem{rem}[thm]{Remark}
 \numberwithin{equation}{section}

 \newcommand{\spec}{{\text{spec }}}
 \newcommand{\pspec}{\text{pspec }}
 \newcommand{\symp}{\text{symp }}
 \newcommand{\prim}{\text{prim }}

\begin{document}

\title[multi-parameter symplectic and Euclidean spaces]{Poisson structures of multi-parameter
symplectic and Euclidean spaces}

\author{Sei-Qwon Oh}

\address{Department of Mathematics, Chungnam National  University,
Taejon 305-764, Korea}

\email{sqoh@math.cnu.ac.kr}

\thanks{The author is supported by the Korea Research Foundation Grant, KRF-2002-015-CP0010.}


\keywords{Poisson algebra, quantized algebra}



\begin{abstract}
A class of Poisson algebras considered as a Poisson version of the multiparameter quantized coordinate rings of
symplectic and Euclidean $2n$-spaces is constructed and
the prime Poisson ideals and the symplectic ideals of these Poisson
algebras are described. As a result,
it is shown that the multiparameter quantized symplectic and Euclidean $2n$-spaces are topological quotients
of their classical spaces.
\end{abstract}

\maketitle

\section*{Introduction}

There is the following conjecture in \cite[II.10.12]{BrGo}: Primitive spectra of quantized algebras
are topological quotients of their classical spaces. This conjecture is known to hold in the cases
$\mathcal{O}_q(SL_2)$, $\mathcal{O}_q(({\bf k}^\times)^n)$ in \cite{Van}, $\mathcal{O}_q({\bf k}^n)$ in \cite{GoLet2}
and $\mathcal{O}_q(\frak{sp}({\bf k}^4))$ in \cite{Hor2}. The main purpose of this paper
is to show that this conjecture is true for the multiparameter quantized symplectic and Euclidean $2n$-spaces.

A quantization of a Poisson algebra is a certain associative algebra with multiplication
deformed by a given Poisson bracket \cite[Ch. 6]{ChPr} and thus quantized spaces seem to be
naturally related to their Poisson structures. For instance, Hodges, Levasseur and Toro described in
\cite{HoLet}
that the primitive ideals of a multiparameter quantum group correspond to the symplectic leaves
of its Poisson variety in the case when they are algebraic and Vancliff have a similar result for
$\mathcal{O}_q(M_2)$ in\cite{Van}.
Goodearl and Letzter proved in \cite{GoLet2}
that the prime and primitive spectra of $\mathcal{O}_{\bf q}({\bf k}^n)$,
the multiparameter quantized coordinate ring of affine $n$-space over an algebraically closed field
are topological quotients of the corresponding classical spectra and the author showed in \cite{Oh4}
that the prime and primitive spectra of $\mathcal{O}_{\bf q}({\bf k}^n)$ are topological quotients of
the corresponding Poisson spectra. Hence it seems that if $A$ is a Poisson algebra
which is the coordinate ring of an affine variety $V$ then the prime and primitive spectra
of standard quantized coordinate rings of $V$ are topological quotients of the prime Poisson and symplectic
spectra of $A$. Here we investigate the Poisson structures for
the multiparameter quantized symplectic and Euclidean $2n$-spaces and then we prove using the Poisson structures
that the prime and primitive spectra of the multiparameter quantized symplectic and Euclidean $2n$-spaces
are topological quotients of their classical corresponding  spaces.

A class of algebras  $K_{n,\Gamma}^{P,Q}$, constructed by Horton in \cite{Hor1}, includes
the multiparameter quantized coordinate rings of symplectic and Euclidean $2n$-spaces,
the graded quantized Weyl algebra, the quantized Heisenberg space,
  and is similar to a class of iterated skew
polynomial rings constructed by G\'omez-Torrecillas and Kaoutit in \cite{GoToKa}. The prime and primitive
spectra for the multiparameter quantized coordinate rings of symplectic and Euclidean $2n$-spaces
were established by G\'omez-Torrecillas and Kaoutit in \cite{GoToKa}, by Horton in \cite{Hor1}
and by the author in \cite{Oh2}.
Here we construct a class of Poisson algebras $A_{n,\Gamma}^{P,Q}$, which is considered as a Poisson version
of $K_{n,\Gamma}^{P,Q}$ and investigate the Poisson structure of $A_{n,\Gamma}^{P,Q}$.
The Poisson structures for $A_{n,\Gamma}^{P,Q}$ obtained here may be considered as Poisson versions
for the algebraic structures of the
multiparameter quantized coordinate rings of symplectic and Euclidean $2n$-spaces
established by G\'omez-Torrecillas and Kaoutit in \cite{GoToKa} or by Horton in \cite{Hor1}.

In the section 1 and 2, we construct a Poisson polynomial ring which is considered as a
Poisson version of a skew polynomial ring and study several basic properties for
the Poisson algebras $A_{n,\Gamma}^{P,Q}$.
In the section 3, we consider an additive group $K$ acting by Poisson derivations on $A_{n,\Gamma}^{P,Q}$
which gives a classification of $K$-prime Poisson ideals of $A_{n,\Gamma}^{P,Q}$. Here we see that
the additive group $K$ is considered as a Poisson version of a multiplicative
group acting by automorphisms on $K_{n,\Gamma}^{P,Q}$.
In the section 4,
we prove that the prime  and primitive spectra of $K_{n,\Gamma}^{P,Q}$ are  topological quotients
of the prime Poisson and symplectic spectra of $A_{n,\Gamma}^{P,Q}$ and, as a corollary, we have
that  the conjecture \cite[II.10.12]{BrGo} for the quantum symplectic and Euclidean $2n$-spaces is true, that is,
 the prime  and primitive spectra of $K_{n,\Gamma}^{P,Q}$ are topological
quotients of the corresponding classical spectra.

\medskip
Assume throughout the paper that ${\bf k}$ denotes an algebraically
closed field of characteristic zero and that all vector spaces are over ${\bf k}$.
 A Poisson algebra $A$ is always a commutative ${\bf k}$-algebra
with ${\bf k}$-bilinear map $\{\cdot,\cdot\}$, called a Poisson bracket,
such that $(A, \{\cdot,\cdot\})$ is a Lie algebra and $\{\cdot,\cdot\}$ satisfies the Leibniz rule, that is,
$$\{ab,c\}=a\{b,c\}+b\{a,c\}$$
for all $a,b,c\in A$. Hence, for any element $a\in A$,
the map
$$h_a:A\longrightarrow A,\ \ h_a(b)=\{a,b\}$$
is a derivation in $A$ which is called a  Hamiltonian defined by $a$.
Assume throughout the paper that $\mathcal{H}(A)$ denotes the set of all Hamiltonians of $A$.

\section{Poisson polynomial ring}

Let $A$ be a Poisson algebra.  A derivation $\delta$ on $A$ is said to be a
 Poisson derivation if
$\delta(\{a,b\})=\{\delta(a),b\}+\{a,\delta(b)\}$ for all $a,b\in
A$.

\begin{thm}
 For a Poisson algebra $A$ with Poisson bracket $\{\cdot,\cdot\}_A$ and
 $\bf k$-linear maps $\alpha,\delta$ from $A$
into itself,
the polynomial ring $A[x]$ is a Poisson algebra
with Poisson bracket
\begin{equation}
\{a,x\}=\alpha(a)x+\delta(a)
\end{equation}
for all $a\in A$ if and only if $\alpha$ is a Poisson derivation and $\delta$ is a derivation such that
\begin{equation}
\delta(\{a,b\}_A)-\{\delta(a),b\}_A-\{a,\delta(b)\}_A=\delta(a)\alpha(b)-\alpha(a)\delta(b)
\end{equation}
for all $a,b\in A$. In this case,  we denote the Poisson algebra $A[x]$  by
$A[x;\alpha,\delta]_p$ and if $\delta=0$ then we simply write
$A[x;\alpha]_p$ for $A[x;\alpha,0]_p$.
\end{thm}
\begin{proof}
If $A[x]$ is a Poisson algebra with the Poisson bracket (1.1) then we have that
$$\aligned\{ab,x\}&=\alpha(ab)x+\delta(ab)\\
a\{b,x\}+\{a,x\}b&=(a\alpha(b)+\alpha(a)b)x+(a\delta(b)+\delta(a)b)\endaligned
$$ for all $a,b\in A$, and thus both $\alpha$ and $\delta$ are derivations on $A$.
Moreover, since the Poisson bracket $\{\cdot,\cdot\}$ satisfies the Jacobi identity, we
have that
$$\aligned
0&=\{\{a,b\},x\}+\{\{b,x\},a\}+\{\{x,a\},b\}\\
&=(\alpha(\{a,b\}_A)-\{\alpha(a),b\}_A-\{a,\alpha(b)\}_A)x\\
&\qquad+\delta(\{a,b\}_A)-\{\delta(a),b\}_A-\{a,\delta(b)\}_A-\delta(a)\alpha(b)+\alpha(a)\delta(b)
\endaligned$$
for all $a,b\in A$. Hence $\alpha$ is a Poisson derivation
and $\delta$ is a derivation such that the pair $(\alpha,\delta)$ satisfies (1.2).

Conversely, we suppose that $\alpha$ is a Poisson derivation and $\delta$ is a derivation
  satisfying the condition (1.2). It is enough to check that
the $\bf k$-bilinear map $\{\cdot,\cdot\}:A[x]\times
A[x]\longrightarrow A[x]$ defined by
\begin{equation}
\aligned\{ax^i,bx^j\}=(\{a,b\}_A&+jb\alpha(a)-ia\alpha(b))x^{i+j}\\
&\qquad+(jb\delta(a)-ia\delta(b))x^{i+j-1}\endaligned
\end{equation} for all monomials $ax^i$ and $bx^j$ in $A[x]$
 is a Poisson bracket on $A[x]$ since
 (1.1) is the case for $i=0, j=1$ and $b=1$ in (1.3)
 and every derivation is uniquely determined by images for generators.

It is easy to check that $\{f,g\}=-\{g,f\}$ for all $f,g\in A[x]$
and that, for a fixed element $f\in A[x]$, the $\bf k$-linear maps
$$ \{f,\cdot\}:A[x]\longrightarrow A[x],\ g\mapsto\{f,g\}{\text{
and }} \{\cdot, f\}:A[x]\longrightarrow A[x],\ g\mapsto\{g,f\}$$
are derivations on $A[x]$. It remains to check that the bracket
$\{\cdot,\cdot\}$ given in (1.3) satisfies the Jacobi identity.  For
$ax^i,bx^j,cx^k\in A[x]$, it is checked that $$\{\{ax^i,bx^j\},
cx^k\}+\{\{bx^j,cx^k\},ax^i\}+\{\{cx^k,ax^i\},bx^j\}=0$$ by using
Leibniz rule and the induction on $i,j,k$. It completes the proof.
\end{proof}

\begin{lem}
Let $\alpha$ and $\delta$ be derivations on a Poisson algebra $A$.
\begin{itemize}
\item[(a)] If $\alpha(a)=\delta(a)$ for all generators $a$ of the
algebra $A$ then $\alpha=\delta$.
\item[(b)] If $\alpha\delta(a)=\delta\alpha(a)$ for all generators
$a$ of the algebra $A$ then $\alpha\delta=\delta\alpha$.
\item[(c)] If $\alpha$ satisfies
the condition $\alpha(\{a,b\})=\{\alpha(a),b\}+\{a,\alpha(b)\}$
for all generators $a,b$ of the algebra $A$ then $\alpha$ is a
Poisson derivation.
\item[(d)] If $\alpha$ and $\delta$ satisfy the condition (1.2) for all
generators of the algebra $A$, then $\alpha$ and $\delta$ satisfy
(1.2) for all elements in $A$.
\end{itemize}
\end{lem}
\begin{proof} We proceed  by induction on the length of monomials of generators in $A$.
Now, it is easy to check (a) and (b) and thus we prove (c) and (d).
Let $a,b,c$ be monomials of generators of $A$.

(c) By the induction hypothesis, we have that
$$\aligned\alpha(\{ab,c\})&=\alpha(a\{b,c\}+b\{a,c\})\\
&=\alpha(a)\{b,c\}+a\alpha(\{b,c\})+\alpha(b)\{a,c\}+b\alpha(\{a,c\})\\
&=\alpha(a)\{b,c\}+a(\{\alpha(b),c\}+\{b,\alpha(c)\})\\
&\qquad+\alpha(b)\{a,c\}+b(\{\alpha(a),c\}+\{a,\alpha(c)\})\\
\{\alpha(ab),c\}+\{ab,\alpha(c)\}&=\{\alpha(a)b,c\}+\{a\alpha(b),c\}+\{ab,\alpha(c)\}\\
&=b\{\alpha(a),c\}+\alpha(a)\{b,c\}+a\{\alpha(b),c\}\\
&\qquad+\alpha(b)\{b,c\}+a\{b,\alpha(c)\}+b\{a,\alpha(c)\}.
\endaligned$$
Hence $\alpha(\{ab,c\})=\{\alpha(ab),c\}+\{ab,\alpha(c)\}$ and so
$\alpha$ is a Poisson derivation by the induction on the length of
monomials.

(d) By the induction hypothesis, we have that
$$\aligned\delta(\{ab,c\})&-\{\delta(ab),c\}-\{ab,\delta(c)\}\\
&=\delta(a)\{b,c\}+a\delta(\{b,c\})+\delta(b)\{a,c\}+b\delta(\{a,c\})\\
&\qquad-\delta(a)\{b,c\}-b\{\delta(a),c\}-\delta(b)\{a,c\}-a\{\delta(b),c\}\\
&\qquad-a\{b,\delta(c)\}-b\{a,\delta(c)\}\\
&=a(\delta(\{b,c\})-\{\delta(b),c\}-\{b,\delta(c)\})\\
&\qquad+b(\delta(\{a,c\})-\{\delta(a),c\}-\{a,\delta(c)\})\\
&=a(\delta(b)\alpha(c)-\alpha(b)\delta(c))+b(\delta(a)\alpha(c)-\alpha(a)\delta(c))\\
\delta(ab)\alpha(c)&-\alpha(ab)\delta(c)\\
&=\delta(a)b\alpha(c)+a\delta(b)\alpha(c)-\alpha(a)b\delta(c)-a\alpha(b)\delta(c).
\endaligned$$
Hence we have
$\delta(\{ab,c\})-\{\delta(ab),c\}-\{ab,\delta(c)\}=\delta(ab)\alpha(c)-\alpha(ab)\delta(c)$,
as claimed.
\end{proof}

An ideal $I$ of a Poisson algebra $A$ is said to be a Poisson ideal of $A$ if $\{I,A\}\subseteq I$.
For Poisson algebras $A$ and $B$, an algebra homomorphism $f:A\longrightarrow B$ is said to be a
Poisson homomorphism if $f(\{a,b\})=\{f(a),f(b)\}$ for all $a,b\in A$.

\begin{lem}
For a Poisson algebra $A$,  let $\alpha$ be a Poisson derivation
and let $\delta$ be a derivation satisfying (1.2).
\begin{itemize} \item[(a)] If $I$ is an
$(\alpha,\delta)$-stable Poisson ideal of $A$ then
$IA[x;\alpha,\delta]_p$ is a Poisson ideal of
$A[x;\alpha,\delta]_p$ and
$A[x;\alpha,\delta]_p/IA[x;\alpha,\delta]_p\cong(A/I)[x;\overline{\alpha},\overline{\delta}]_p$,
where $\overline{\alpha}$ and $\overline{\delta}$ are the maps in $A/I$ induced by
 $\alpha$ and $\delta$ respectively.
 \item[(b)] Let $\beta$ be a Poisson derivation on $A[y;\alpha]_p$ such that
 $\beta(A)\subseteq A$ and $\beta(y)=cy$ for some $c\in{\bf k}$.
 Then $A[y;\alpha]_p[x;\beta]_p\cong
 A[x;\beta|_A]_p[y;\alpha']_p$, where $\alpha'$ is a Poisson derivation on
 $A[x;\beta|_A]_p$ such that $\alpha'|_A=\alpha$ and $\alpha'(x)=-cx$.
 \end{itemize}
\end{lem}
\begin{proof}
(a) Note that the Poisson algebras $A[x;\alpha,\delta]_p$ and
$(A/I)[x;\overline{\alpha},\overline{\delta}]_p$ are constructed
by Theorem 1.1. The map
$\psi:A[x;\alpha,\delta]_p\longrightarrow(A/I)[x;\overline{\alpha},\overline{\delta}]_p$
defined by
$\psi(a_0+a_1x+\cdots+a_nx^n)=\overline{a}_0+\overline{a}_1x+\cdots+\overline{a}_nx^n$
is a Poisson epimorphism and  has the kernel
$IA[x;\alpha,\delta]_p$, hence we have the conclusion.

(b) Observe that there exists a Poisson algebra
$A[x;\beta|_A]_p$ since $\beta|_A$ is a Poisson derivation on $A$.
Since $A[x,y]=A[y;\alpha]_p[x;\beta]_p$ is a Poisson algebra,
the map $\alpha'$ on $A[x;\beta|_A]_p$ satisfying
$\{f,y\}=\alpha'(f)y$ for all $f\in A[x;\beta|_A]_p$ is  a Poisson
derivation  such that $\alpha'|_A=\alpha$ and $\alpha'(x)=-cx$ by
Theorem 1.1, and thus the Poisson algebra
$A[x;\beta|_A]_p[y;\alpha']_p$ is constructed. Clearly, the
identity map from $A[y;\alpha]_p[x;\beta]_p$ into $
A[x;\beta|_A]_p[y;\alpha']_p$ is a Poisson isomorphism.
\end{proof}

\begin{prop}
Let $A$ be a Poisson algebra. For Poisson derivations $\alpha$ and
$\beta$ on $A$, $c\in {\bf k}$ and $u\in A$ such that
$$\alpha\beta=\beta\alpha, \ \ \{a,u\}=(\alpha+\beta)(a)u$$
for all $a\in A$, the polynomial ring $A[y,x]$ has the following
Poisson bracket
\begin{equation}
\{a,y\}=\alpha(a)y,\  \{a,x\}=\beta(a)x,\  \{y,x\}=cyx+u
\end{equation}
 for all $a\in A$. The Poisson algebra $A[y,x]$ with Poisson bracket (1.4)
 can be presented by $A[y;\alpha]_p[x;\beta',\delta]_p$, where
 $\beta'$ is the Poisson derivation on $A[y;\alpha]_p$ such that
$\beta'|_A=\beta$ and $\beta'(y)=cy$, and $\delta$ is the derivation on
$A[y;\alpha]_p$ such that $\delta|_A=0$, $\delta(y)=u$.
 We often denote by $(A;\alpha,\beta, c,u)$ the Poisson algebra
$A[y,x]$ with Poisson bracket (1.4).
\end{prop}
\begin{proof}
By Theorem 1.1, there exists the Poisson algebra $A[y;\alpha]_p$
with Poisson bracket $\{a,y\}=\alpha(a)y$ for all $a\in A$ and the
derivation $\beta$ is extended to a derivation,  denoted by
$\beta'$, to $A[y;\alpha]_p$ by setting $\beta'(y)=cy$. Note that
the derivation $\delta=u\frac{d}{dy}$ on $A[y;\alpha]_p$ satisfies
$\delta(y)=u$ and $\delta(a)=0$ for all $a\in A$. Let us prove
that, for all $f,g\in A[y;\alpha]_p$,
\begin{equation}\aligned \beta'(\{f,g\})&=\{\beta'(f),g\}+\{f,\beta'(g)\} \\
\delta(\{f,g\})&=\{\delta(f),g\}+\{f,\delta(g)\}+\delta(f)\beta'(g)-\beta'(f)\delta(g).
\endaligned
\end{equation}

If $f,g\in A$ then the formulas in (1.5) hold trivially since $\beta'$ is a
Poisson derivation on $A$. Hence it is
enough to prove (1.5) for the case $f=a\in A$ and $g=y$ by
Lemma 1.2. Now we have that $$\aligned
\beta'(\{a,y\})&=\beta'(\alpha(a)y)=\alpha(a)\beta'(y)+\beta'(\alpha(a))y
\\ &=c\alpha(a)y+\alpha(\beta(a))y=\{\beta'(a),y\}+\{a,\beta'(y)\}\\
\delta(\{a,y\})&=\delta(\alpha(a)y)=\alpha(a)u=\{a,u\}-\beta(a)u\\
&=\{\delta(a),y\}+\{a,\delta(y)\}+\delta(a)\beta'(y)-\beta'(a)\delta(y),\endaligned$$
as claimed.

Therefore $\beta'$ is a Poisson derivation on $A[y;\alpha]_p$ such
that the pair $(\beta',\delta)$ satisfies (1.2), and thus,
by Theorem 1.1, there exists the Poisson algebra
$A[y,x]=A[y;\alpha]_p[x;\beta',\delta]_p$ with the Poisson bracket (1.4).
\end{proof}

\medskip
\noindent {\bf Example.} For the Poisson algebra ${\bf k}[b,c]$
with trivial Poisson bracket, that is, $\{b,c\}=0$, the derivation
$\alpha=-2b\frac{\partial}{\partial b}-2c\frac{\partial}{\partial c}$ on
${\bf k}[b,c]$ is a Poisson derivation clearly. Observe that
 the Poisson algebra $({\bf
k}[b,c],\alpha,-\alpha,0, 4bc)$ is the Poisson algebra ${\bf
k}[b,c][a,d]$ with Poisson bracket
\[\begin{tabular}{lll}
$\{b,c\}=0$, &$\{b,a\}=-2ba$, &$ \{c,a\}=-2ca$,\\
$\{b,d\}=2bd$, &$\{c,d\}=2cd$, &$\{a,d\}=4bc$,\end{tabular}\]
 which is the
Poisson algebra given in \cite[2.9]{Oh4}, \cite[Example
3.2.9]{KoSo} and \cite[3.13]{Van}.

\medskip
\noindent {\bf Example.} The Poisson algebra $({\bf k}, 0,0,0,1)$
is the algebra ${\bf k}[y,x]$ with Poisson bracket $\{y,x\}=1$,
which is equal to the Poisson algebra ${\bf k}[y,x]$ with Poisson
bracket $\{f,g\}=\frac{\partial f}{\partial y}\frac{\partial
g}{\partial x}- \frac{\partial g}{\partial y}\frac{\partial
f}{\partial x}$, $f,g\in {\bf k}[y,x]$, given in \cite[p. 18]
{ChPr}.

\begin{lem}
Let $A$ and $B$ be commutative ${\bf k}$-algebras, $S$ a multiplicative set of $A$ and
$\psi:A\longrightarrow B$ an algebra homomorphism.
 If $\delta$ is a ${\bf k}$-linear map
from $A$ into $B$ such that \begin{equation}\delta(ab)=
\delta(a)\psi(b)+\psi(a)\delta(b)\end{equation} for all $a,b\in
A$, then there exists a unique ${\bf k}$-linear map $\delta'$ from
$A[S^{-1}]$ into $B[\psi(S)^{-1}]$ such that
$\delta'(a)=\delta(a)$ for $a\in A$ and
$$\delta'((as^{-1})(bt^{-1}))=\delta'(as^{-1})\psi'(bt^{-1})+\psi'(as^{-1})\delta(bt^{-1})$$
for all $as^{-1}, bt^{-1}\in A[S^{-1}]$, where $\psi'$ is the
extension of $\psi$ from $A[S^{-1}]$ into $B[\psi(S)^{-1}]$.
\end{lem}
\begin{proof}
Define
$\delta'(as^{-1})=(\delta(a)\psi(s)-\psi(a)\delta(s))\psi(s)^{-2}$
for all $as^{-1}\in A[S^{-1}]$.  It is checked routinely that
$\delta'$ is a well-defined ${\bf k}$-linear map and satisfies the
required conditions.
\end{proof}

\begin{prop}
Let $S$ be a multiplicative set of a Poisson algebra $A$.
\begin{itemize}
 \item[(a)] The localization
$A[S^{-1}]$ has the Poisson bracket defined by
$$\{as^{-1},bt^{-1}\}=(\{a,b\}st-\{a,t\}bs-\{s,b\}at+\{s,t\}ab)s^{-2}t^{-2}$$
for $as^{-1},bt^{-1}\in  A[S^{-1}]$.
\item[(b)] If $\delta$ is a (respectively,
Poisson) derivation on $A$ then there exists a (respectively,
Poisson) derivation $\delta'$ on $S^{-1}A$ defined by
 $\delta'(as^{-1})=(\delta(a)s-a\delta(s))s^{-2}$
 for $as^{-1}\in A[S^{-1}]$.
 \item[(c)] If a group $H$ acts by (respectively, Poisson)
 derivations on $A$ then $H$ acts by (respectively, Poisson)
 derivations on $A[S^{-1}]$. Moreover, if $a\in A$ and $s\in S$ are
 $H$-eigenvectors and $as^{-1}\neq0$ then $as^{-1}\in A[S^{-1}]$ is also an
 $H$-eigenvector.
\end{itemize}
\end{prop}
\begin{proof}
 If $\delta$ is a derivation on $A$ then, by \cite[14.2.2]{McRo}, there
exists a unique derivation $\delta'$ on $S^{-1}A$ such that
$\delta'(a)=\delta(a)$ for all $a\in A$, and thus
$\delta'(as^{-1})=(\delta(a)s-a\delta(s))s^{-2}$
 for $as^{-1}\in A[S^{-1}]$.

(a) For $a\in A$, let $\ell_a$ be the derivation on $A$ defined by
$\ell_a(b)=\{a,b\}$ for all $b\in A$. Then there exists a
derivation $\ell'_a$ on $A[S^{-1}]$ such that
$\ell'(bs^{-1})=(\ell_a(b)s-b\ell_a(s))s^{-2}$
 for $bs^{-1}\in A[S^{-1}]$. Fix an element $bs^{-1}\in A[S^{-1}]$
 and let $j$ be the natural homomorphism from $A$ into
 $A[S^{-1}]$. Then the ${\bf k}$-linear map $r_{bs^{-1}}$ from $A$ into $A[S^{-1}]$ defined by
 $r_{bs^{-1}}(a)=\ell_a(bs^{-1})$ for all $a\in A$ satisfies (1.6), and thus there exists  a derivation
 $r'_{bs^{-1}}$ on $A[S^{-1}]$ such that $r'_{bs^{-1}}(a)=r_{bs^{-1}}(a)$  by Lemma 1.5 since
 $A[S^{-1}][j(S)^{-1}]=A[S^{-1}]$.

Now, define the Poisson bracket on $A[S^{-1}]$ by
 $$\{as^{-1}, bt^{-1}\}=r'_{bt^{-1}}(as^{-1})$$
 for $as^{-1},bt^{-1}\in A[S^{-1}]$. Then it is easy to check that
 $A[S^{-1}]$ is a Poisson algebra and that $$\{as^{-1},bt^{-1}\}=(\{a,b\}st-\{a,t\}bs-\{s,b\}at+\{s,t\}ab)s^{-2}t^{-2}$$
for $as^{-1},bt^{-1}\in  A[S^{-1}]$.

(b) By \cite[14.2.2]{McRo}, $\delta'$ is the unique derivation such that $\delta'(a)=\delta(a)$ for all $a\in A$
and if  $\delta$ is a Poisson derivation  then the fact that $\delta'$ is a also a Poisson derivation is
verified immediately by using Lemma 1.3.

(c) It follows immediately from (b).
\end{proof}

An element
$z$ of a Poisson algebra $A$ is said to be Poisson normal if $z$ is not a zero divisor and
$\{z,A\}\subseteq zA$. Hence
if $z\in A$ is Poisson normal then there exists a Poisson derivation $\gamma$ on $A$ defined by
$$\{a,z\}=\gamma(a)z$$
for all $a\in A$.

\begin{lem}
In the Poisson algebra $B=(A;\alpha,\beta,c,u)$ given in Proposition 1.4, suppose that
$\alpha(u)=du,\beta(u)=-du$ for some $d\in{\bf k}$ such that $c+d\neq0$ and set $$z=(c+d)yx+u.$$ Then
\begin{itemize}
\item[(a)]
$z$ is a Poisson normal element of $B$. More precisely,
$$\{a,z\}=(\alpha+\beta)(a)z, \ \ \{y,z\}=cyz, \ \ \{x,z\}=-cxz$$
for all $a\in A$.
\item[(b)]
$B[y^{-1}]=A[y^{\pm1};\alpha]_p[y^{-1}z;\beta']_p$, where $\beta'$ is the extension of $\beta$
with $\beta'(y)=cy$.
\end{itemize}
\end{lem}
\begin{proof}
(a) It follows immediately from a straight  calculation.

(b) Note that $B[y^{-1}]=A[y^{\pm1};\alpha]_p[x;\beta',\delta]_p$,
where $\beta'$ and $\delta$ are given in
 Proposition 1.4.
Since
$$\aligned
\{a,y^{-1}z\}&=\beta(a)y^{-1}z, \ \ \ a\in A\\
\{y,y^{-1}z\}&=cy(y^{-1}z)\\
\{x,y^{-1}z\}&=y^{-1}u(y^{-1}z),
\endaligned$$
 $y^{-1}z$ is a Poisson normal element of
 $B[y^{-1}]$ and  $\{f, y^{-1}z\}=\beta'(f)y^{-1}z$ for all $f\in A[y^{\pm1};\alpha]_p$.
Hence the Poisson algebra $A[y^{\pm1};\alpha]_p[y^{-1}z;\beta']_p$ is constructed by Theorem 1.1 and
 $B[y^{-1}]=A[y^{\pm1};\alpha]_p[y^{-1}z;\beta']_p$ since
$$x=(c+d)^{-1}y^{-1}z-(c+d)^{-1}y^{-1}u\in A[y^{\pm1};\alpha]_p[y^{-1}z;\beta']_p.$$
\end{proof}

Let $A$ be a Poisson algebra. A prime ideal of $A$ which is also a Poisson ideal is called a prime Poisson ideal
and a Poisson ideal $N$ of $A$ is said to be symplectic if there exists a maximal
ideal $M$ such that $N$ is the largest Poisson ideal contained in $M$ (\cite[Definition 1.2]{Oh4}). For
an ideal $I$ of  $A$, denote by $(I:\mathcal{H}(A))$ the largest $\mathcal{H}(A)$-stable ideal
contained in $I$. Note that $(I:\mathcal{H}(A))$ is the largest Poisson ideal contained in $I$
since $\mathcal{H}(A)$ is the set of all Hamiltonians of $A$. Hence if $M$ is a maximal ideal of $A$ then
$(M:\mathcal{H}(A))$ is a symplectic ideal.

\begin{lem}
Let $A$ be a Poisson algebra and let $I$ be a proper Poisson ideal of $A$. Then there exists
a prime Poisson ideal $P$ of $A$ containing $I$.
\end{lem}
\begin{proof}
Since $I$ is proper, there exists a prime ideal $Q$ of $A$ such that $I\subseteq Q$.
Set $P=(Q:\mathcal{H}(A))$.
Then $P$ is a prime Poisson ideal containing $I$ by \cite[Lemma 1.3]{Oh4}.
\end{proof}

\begin{defn}
For a ${\bf k}$-algebra $R$ and a Poisson ${\bf k}$-algebra $A$, set
$$\aligned
\text{spec}(R)&=\text{the set of all prime ideals of $R$} \\
\text{prim}(R)&=\text{the set of all primitive ideals of $R$}\\
\text{pspec}(A)&=\text{the set of all prime Poisson ideals of $A$}\\
\text{symp}(A)&=\text{the set of all symplectic ideals of $A$}\\
\text{max}(A)&=\text{the set of all maximal ideals of $A$}.
\endaligned$$
The sets $\text{spec}(R)$ and $\text{spec}(A)$ are topological spaces equipped with Zariski topologies
and the others are also topological spaces equipped
with relative topologies since $\text{prim}(R)\subseteq \text{spec}(R)$, $\max A\subseteq\spec A$ and
$\text{symp}(A)\subseteq\text{pspec}(A)\subseteq \text{spec}(A)$.
\end{defn}

\section{Poisson algebra $A_n=A_{n,\Gamma}^{P,Q}$}

\begin{thm} Let $\Gamma=(\gamma_{ij})$ be a skew-symmetric $n\times n$-matrix
with entries in ${\bf k}$, that is, $\gamma_{ij}=-\gamma_{ji}$ for all $i, j=1,\cdots,n$. Let
$P=(p_1,p_2,\cdots,p_n)$ and $Q=(q_1,q_2,\cdots,q_n)$ be elements
of ${\bf k}^{n}$ such that $p_i\neq q_i$ for each $i=1,\cdots,n$.
Then the polynomial ring ${\bf k}[y_1,x_1,\cdots,y_n,x_n]$ has the
following Poisson bracket:
\begin{equation}\begin{array}{ll}
\{y_i,y_j\}=\gamma_{ij}y_iy_j &(\mbox{all } i,j)\\
\{x_i,y_j\}=(p_j-\gamma_{ij})y_jx_i  &(i<j)\\
\{y_i,x_j\}=-(q_i+\gamma_{ij})y_ix_j &(i<j)\\
\{x_i,x_j\}=(q_i-p_j+\gamma_{ij})x_ix_j  &(i<j)\\
\{x_i,y_i\}=q_iy_ix_i+\sum_{k=1}^{i-1}(q_k-p_k)y_kx_k &(\mbox{all
} i)
 \end{array}
 \end{equation}
The Poisson algebra ${\bf k}[y_1,x_1,\cdots,y_n,x_n]$ is denoted
by $A_{n,\Gamma}^{P,Q}$ or by $A_n$ unless any confusion arises.
\end{thm}
\begin{proof}
Let ${\bf k}[h]$ be the polynomial ring and let $A$ be the ${\bf
k}[h]$-algebra generated by $y'_1,x'_1,\cdots, y'_n,x'_n$ subject
to the following relations:
\[\begin{array}{ll}
y'_iy'_j-y'_jy'_i=\gamma_{ij}hy'_jy'_i &(\mbox{all } i,j)\\
x'_iy'_j-y'_jx'_i=(p_j-\gamma_{ij})hy'_jx'_i  &(i<j)\\
y'_ix'_j-x'_jy'_i=-(q_i+\gamma_{ij})hx'_jy'_i &(i<j)\\
x'_ix'_j-x'_jx'_i=(q_i-p_j+\gamma_{ij})hx'_jx'_i  &(i<j)\\
x'_iy'_i-y'_ix'_i=q_ihy'_ix'_i+\sum_{k=1}^{i-1}(q_k-p_k)hy'_kx'_k
&(\mbox{all } i)
 \end{array}\]
 Then $h$ is not a zero divisor of $A$ and $A_n=A/hA$ is the commutative polynomial ring ${\bf
 k}[y_1,x_1,\cdots,y_n,x_n]$, where
 $$y_i=y'_i+hA,\ \ \  x_i=x'_i+hA$$
 for all $i=1,2,\cdots,n$.
 Moreover, the algebra $A_n$ has the Poisson bracket (2.1)
 by \cite[III.5.4]{BrGo}. It completes the proof.
\end{proof}

\begin{rem} Set
$$A_0={\bf k},\ \ \ \ \ A_j={\bf k}[y_1,x_1,\cdots, y_j,x_j]\subseteq A_{n,\Gamma}^{P,Q}$$
for
each $j=0,1,\cdots, n$. Then each $A_j$ is a Poisson subalgebra of
$A_{n,\Gamma}^{P,Q}$ and $A_j=A_{j-1}[y_j,x_j]$ for each $j$, and thus, by
Theorem 1.1, there exist Poisson derivations
$\alpha_{j},\beta_{j}$ and a derivation $\delta_j$
such that $A_j$ can be presented by
$$A_j=A_{j-1}[y_j;\alpha_j]_p[x_j;\beta_j,\delta_j]_p,$$ where
\begin{equation}\begin{array}{lll}
\alpha_{j}(y_i)=\gamma_{ij}y_i, &\alpha_{j}(x_i)=(p_j-\gamma_{ij})x_i &(i<j)\\
 \beta_{j}(y_i)=-(q_i+\gamma_{ij})y_i,
 &\beta_{j}(x_i)=(q_i-p_j+\gamma_{ij})x_i&(i<j)\\
 \delta_j(y_i)=0,&\delta_j(x_i)=0&(i<j)\\
 \beta_j(y_{j})=-q_{j}y_{j}&\delta_j(y_{j})=-\sum_{k=1}^{j-1}(q_k-p_k)y_kx_k.&
 \end{array}
 \end{equation}
 Set
  $$\Omega_0=0,\ \ \ \ \
 \Omega_{j}=\sum_{k=1}^{j}(q_k-p_k)y_kx_k$$
for all $j=1,\cdots,n-1$, and note that
$$\alpha_j\beta_j=\beta_j\alpha_j,   \ \ \
\{a,\Omega_{j-1}\}=(\alpha_j+\beta_j)(a)\Omega_{j-1}$$ for all $a\in A_{j-1}$.
 Hence we have $A_j=(A_{j-1};\alpha_j,\beta_j, -q_j, -\Omega_{j-1})$ by Lemma 1.4
  and so the Poisson algebra
$A_n=A_{n,\Gamma}^{P,Q}$ has the chain of Poisson subalgebras
$$A_0={\bf k}\subseteq A_1=(A_0;\alpha_1,\beta_1,-q_1,0)\subseteq
\cdots\subseteq A_n=(A_{n-1};\alpha_n,\beta_n,-q_n,-\Omega_{n-1}).$$
\end{rem}

\begin{lem}
As in Remark 2.2, set
$$\Omega_{i}=\sum_{k=1}^{i}(q_k-p_k)y_kx_k\in A_n=A_{n,\Gamma}^{P,Q}$$ for each $i=1,\cdots
,n$ and $\Omega_0=0$.

(a) For any $\Omega_j$,
\[\begin{array}{lll} \{y_i,\Omega_j\}=-q_iy_i\Omega_j, &\{x_i,\Omega_j\}=q_ix_i\Omega_j, &(i\leq
j)\\
\{y_i,\Omega_j\}=-p_iy_i\Omega_j, &\{x_i,\Omega_j\}=p_ix_i\Omega_j, &(i>j)\\
\{\Omega_i,\Omega_j\}=0, & &({\mbox{all }} i,j)
\end{array}\]

(b) We have the following relations:
\begin{equation}
\Omega_{i-1}=\{x_i,y_i\}-q_iy_ix_i, \ \ \
\Omega_i=\{x_i,y_i\}-p_iy_ix_i
\end{equation}
 Hence, $y_i$ and $x_i$ are
Poisson normal modulo $\langle \Omega_i\rangle$ and $\langle
\Omega_{i-1}\rangle$.
\end{lem}
\begin{proof}
The formulas of (a) follow from (2.1) and the
formulas of (b) follow immediately since
$\Omega_i=(q_i-p_i)y_ix_i+\Omega_{i-1}$ and
$\{x_i,y_i\}=q_iy_ix_i+\Omega_{i-1}$.
\end{proof}

\begin{defn} \cite[Definition 1.4]{Oh2} Let
$\mathcal P_n=\{\Omega_1, y_1,x_1,\cdots,
\Omega_n,y_n,x_n\}\subseteq A_{n}.$ A subset $T$ of
$\mathcal P_n$ is said to be {\it admissible} if it satisfies the
conditions:
\begin{itemize}
\item[(a)] $y_i$ or $x_i\in T$ $\Leftrightarrow$ $\Omega_i$ and
$\Omega_{i-1}\in T$   \ \ \ $(2\leq i\leq n)$
\item[(b)] $y_1$ or
$x_1\in T$ $\Leftrightarrow$ $\Omega_1\in T$.
\end{itemize}
\end{defn}

\begin{defn}
We define an order relation on the generators of $A_n$ by
$$y_1,\  x_1,\  y_2, \ x_2, \ \cdots,\  y_n,\  x_n.$$
Hence the standard monomials of $A_n$ are of the form
$y_1^{r_1}x_1^{s_1}y_2^{r_2}x_1^{s_2}\cdots y_n^{r_n}x_n^{s_n}$,
where $r_i,s_i$ are nonnegative integers.
\end{defn}

Let $T$ be an admissible set. In order to find  a ${\bf k}$-basis for $A_n/\langle T\rangle$,
we use an argument for a Gr\"obner-Shirshov  basis.
Refer to \cite{OhPaSh2} and \cite{KaLe1} for  further background and terminologies on
the Gr\"obner-Shirshov  basis.

\begin{lem}
(a) For every admissible set $T$ of $A_n$, $T$ is a Gr\"obner-Shirshov  basis.

(b) The algebra $A_n/\langle T\rangle$ has a ${\bf k}$-basis
consisting of the natural images of all the standard monomials which are not divided by any element
in the set
$$\mathcal A_T=\{y_i\ |\ y_i\in T\}\cup\{x_i\ |\ x_i\in T\}\cup\{y_ix_i\ |\ \Omega_i\in T
, y_i\notin T,x_i\notin T\}.$$
\end{lem}
\begin{proof}
We use the notation given in \cite[1.2]{OhPaSh2}. By \cite[Theorem 1.5]{OhPaSh2}, it is enough to
show that $T$ is closed under composition. Since the maximal monomial of $\Omega_i$ is $y_ix_i$,
a composition occurs for the following three cases:
$$\{y_i,\Omega_i\}\subseteq T,\ \  \{x_i,\Omega_i\}\subseteq T,\ \ \{y_i,x_i,\Omega_i\}\subseteq T.$$
Clearly, all the three cases are closed under composition, hence
 $T$ is a Gr\"obner-Shirshov  basis and (b) follows immediately from \cite[Theorem 1.5]{OhPaSh2}
since all standard monomials form a ${\bf k}$-basis of $A_n$.
\end{proof}

\begin{prop}
For every admissible set $T$, the ideal $\langle T\rangle$ is a
prime Poisson ideal of $A_n$.
\end{prop}
\begin{proof}
We proceed by induction on $n$. If $n=1$ then there are four
admissible sets, namely, $\phi,\{y_1,\Omega_1\},\{x_1,\Omega_1\}$
and $\{y_1,x_1,\Omega_1\}$, which respectively generate the prime
Poisson ideals $0,\langle y_1\rangle,\langle x_1\rangle$ and
$\langle y_1,x_1\rangle$.

Suppose now that $n>1$ and that for each $k<n$, if $U$ is an
admissible set of $A_k$ then the ideal generated by $U$ is a prime
Poisson ideal of $A_k$. Given an admissible set $T$ of $A_n$, let
$T'=T\cap\mathcal P_{n-1}$. Then $T'$ is an admissible set of
$A_{n-1}$ and the ideal $\langle T'\rangle$ is a prime
Poisson ideal of $A_{n-1}$ by the induction hypothesis. Note that
$T$ is one of the following five sets:
$$T',\  T'\cup\{\Omega_n\},\  T'\cup\{y_n,\Omega_n\},\ T'\cup\{x_n,\Omega_n\},\
T'\cup\{y_n,x_n,\Omega_n\}.$$ Since
\begin{equation}\aligned
A_n/\langle T'\rangle
A_n&=A_{n-1}[y_n;\alpha_n]_p[x_n;\beta_n,\delta_n]_p/\langle T'\rangle
A_n\\
&\cong (A_{n-1}/\langle
T'\rangle,\overline\alpha_n,\overline{\beta_n},-q_n,-\overline{
\Omega}_{n-1})=B\endaligned
\end{equation}
by Remark 2.2 and Lemma 1.3, it is enough to prove that the
canonical images of the above five sets in $B$ generate prime
Poisson ideals in $B$.

The case $T=T'$: The ideal of $B$ generated by the canonical image
of $T'$ is 0 clearly.

The case $T=T'\cup\{\Omega_n\}$: Note that $\Omega_{n-1}\notin T'$
and $\Omega_n=(q_n-p_n)y_nx_n+\Omega_{n-1}$. Since
$\Omega_{n-1}\notin \langle T'\rangle$ the canonical image
$\overline \Omega_{n-1}$ in $A_{n-1}/\langle T'\rangle$ is nonzero
and the
 canonical image of $\Omega_n$ in $B$ is $(q_n-p_n)y_nx_n+\overline
 \Omega_{n-1}$, which generates the prime Poisson ideal $\langle y_nx_n+ (q_n-p_n)^{-1}\overline\Omega_{n-1}\rangle$ in $B$.

The case  $T=T'\cup\{y_n,\Omega_n\}$: Since $\Omega_{n-1}\in
 T'$, the canonical image of $\Omega_n$ in $B$ is $(q_n-p_n)y_nx_n$ and
 $\{y_n,x_n\}=-q_ny_nx_n$ in $B$. Hence the ideal of $B$ generated by the
 canonical images of $T$ in $B$ is $\langle y_n\rangle$ which is a
 prime Poisson ideal of $B$.

The case $T=T'\cup\{x_n,\Omega_n\}$: As in the case
 $T=T'\cup\{y_n,\Omega_n\}$, the ideal of $B$ generated by the
 canonical image of $T$ is $\langle x_n\rangle$ which is a
 prime Poisson ideal of $B$.

The case  $T=T'\cup\{y_n,x_n,\Omega_n\}$: Since $\Omega_{n-1}\in
 T'$, the canonical image of $\Omega_n$ in $B$ is $(q_n-p_n)y_nx_n$ and
 $\{y_n,x_n\}=-q_ny_nx_n$ in $B$. Hence the ideal of $B$ generated by the
 canonical images of $T$ in $B$ is $\langle y_n, x_n\rangle$ which is a
 prime Poisson ideal of $B$.
\end{proof}

\begin{prop}
(a) For every prime Poisson ideal $P$ of $A_n$, $P\cap\mathcal
P_n$ is an admissible set.

(b) For an admissible set $T$, let
$$\mbox{pspec}_T
A_n=\{P\in\mbox{pspec}A_n\ |\ P\cap\mathcal P_n=T\}.$$ Then
$\mbox{pspec}A_n$ is the disjoint union of all $\mbox{pspec}_T
A_n$, that is,
$$\mbox{pspec}A_n=\biguplus_{T}\mbox{pspec}_T
A_n.$$
\end{prop}
\begin{proof}
(a) For convenience, set $\Omega_0=0$. If
$\Omega_i,\Omega_{i-1}\in P$ then $y_ix_i\in P$ by (2.3) and so $y_i\in P$
or $x_i\in P$ since $P$ is a prime ideal. Conversely, if $y_i\in
P$ or $x_i\in P$ then $\Omega_i,\Omega_{i-1}\in P$ by (2.3).
Hence we have that, for $i\geq2$, $\Omega_i,\Omega_{i-1}\in
P\cap\mathcal P_n$ if and only if $y_i\in P\cap\mathcal P_n$ or
$x_i\in P\cap\mathcal P_n$, and that $\Omega_1\in P\cap\mathcal P_n$ if
and only if $y_1\in P\cap\mathcal P_n$ or $x_1\in P\cap\mathcal
P_n$. It follows that $P\cap\mathcal P_n$ is admissible.

(b) It follows immediately from (a).
\end{proof}

\begin{defn} \cite[Definition 3.1]{Oh2}
Let $T$ be an admissible set of $A_n$ and let $S_T$ be the subset
of $\mathcal P_n$ consisting of
$$S_T=(\{y_1,x_1,\cdots,y_n,x_n\}\cap T)\cup\{\Omega_i\ |\ \Omega_i\in
T, y_i\notin T, x_i\notin T\}.$$ We define $\text{length}(T)$ to
be the number of elements in $S_T$. Note that the number of elements in $S_T$
is equal to that of $\mathcal A_T$ given in Lemma 2.6.
\end{defn}

\begin{lem} For any admissible set $T$ of $A_n$,
the Gelfand-Kirillov dimension of $A_n/\langle T\rangle$ is equal
to $2n-\text{length}(T)$.
\end{lem}
\begin{proof}
We proceed by induction on $n$. If $n=1$ then there are four admissible sets, namely,
$\phi,\{y_1,\Omega_1\},\{x_1,\Omega_1\}, \{y_1,x_1,\Omega_1\}$. Let $I$ be the ideal of $A_1$ generated by one
of the above sets. Then $A_1/I$ is one of the forms $A_1={\bf k}[y_1,x_1], {\bf k}[x_1],
{\bf k}[y_1], {\bf k}$. Hence the conclusion for $n=1$ is true clearly.

Suppose now that $n>1$ and that for each $k<n$, if $U$ is an admissible set of $A_k$ then
the Gelfand-Kirillov dimension of $A_k/\langle U\rangle$ is equal
to $2k-\text{length}(U)$. Given an admissible set $T$ of $A_n$, let $T'=T\cap \mathcal{P}_{n-1}$. Then
$T'$ is an admissible set of $A_{n-1}$ and the Gelfand-Kirillov dimension of $A_{n-1}/\langle T'\rangle$ is equal
to $2n-2-\text{length}(T')$ by the induction hypothesis.
Note that
$T$ is one of the following five sets:
$$T',\  T'\cup\{\Omega_n\},\  T'\cup\{y_n,\Omega_n\},\ T'\cup\{x_n,\Omega_n\},\
T'\cup\{y_n,x_n,\Omega_n\}.$$
To find the Gelfand-Kirillov dimension of $A_n/\langle T\rangle$, we will use (2.4).
Note that the Gelfand-Kirillov dimension of  $B$ is equal to
 $2n-\text{length}(T')$ by \cite[Example 3.6]{KrLe}.
Let $I$ be the ideal of $B$ generated by the canonical images of the above five sets.

The case $T=T'$: Then $I=0$, $\text{length}(T)=\text{length}(T')$
 and the Gelfand-Kirillov dimension  of $B/I$ is equal to that of $B$, which is equal to
 $2n-\text{length}(T)$.

 The case $T=T'\cup\{\Omega_n\}$: By the proof of Proposition 2.7,
 $I$ is the prime ideal $\langle y_nx_n+ (q_n-p_n)^{-1}\overline\Omega_{n-1}\rangle$ and
 $\text{length}(T)=\text{length}(T')+1$.
 The Gelfand-Kirillov dimension  of $B/I$ is less than or equal to $2n-\text{length}(T')-1$
 by \cite[Corollary 3.16]{KrLe} and larger than or equal to $2n-\text{length}(T')-1$
 since the subalgebra $(A_{n-1}/\langle T'\rangle)[y_n]$ of $B$ has the Gelfand-Kirillov dimension
 $2n-\text{length}(T')-1$
 and the canonical map $(A_{n-1}/\langle T'\rangle)[y_n] \longrightarrow B/I$ is injective.
 Hence the Gelfand-Kirillov dimension of $B/I$ is equal to $2n-\text{length}(T)$.

 The case $T=T'\cup\{y_n,\Omega_n\}$: By the proof of Proposition 2.7, $I=\langle y_n\rangle$ and
 $B/I\cong (A_{n-1}/\langle T'\rangle)[x_n]$ has the Gelfand-Kirillov dimension  $2n-\text{length}(T)$
 since
 $\text{length}(T)=\text{length}(T')+1$.

 The case $T=T'\cup\{x_n,\Omega_n\}$: By the proof of Proposition 2.7, $I=\langle x_n\rangle$ and
 $B/I\cong (A_{n-1}/\langle T'\rangle)[y_n]$ has the Gelfand-Kirillov dimension  $2n-\text{length}(T)$ since
 $\text{length}(T)=\text{length}(T')+1$.

 The case $T=T'\cup\{y_n, x_n,\Omega_n\}$: By the proof of Proposition 2.7, $I=\langle y_n,x_n\rangle$ and
 $B/I\cong (A_{n-1}/\langle T'\rangle)$ has the Gelfand-Kirillov dimension  $2n-\text{length}(T)$
 since
 $\text{length}(T)=\text{length}(T')+2$.

\end{proof}

\section{$K$-actions on $A_{n,\Gamma}^{P,Q}$}

In this section, we will show that every $K$-prime Poisson ideal of
$A_{n,\Gamma}^{P,Q}$ is generated by an admissible set. The statements  and proofs
of this section are modified from those of \cite[\S3]{Hor1}.

\begin{defn} Let
$$\aligned K=\{(h_1,h_2,\cdots,&h_{2n-1},h_{2n})\in {\bf k}^{2n}\ | \\
&\ h_{2i-1}+h_{2i}=h_{2j-1}+h_{2j} {\mbox{ for all }
}i,j=1,\cdots,n\}.\endaligned$$ The additive group $K$ acts on
$A_n$ as follows:
$$(h_1,h_2,\cdots,h_{2n-1},h_{2n})(f)=\sum_i(h_{2i-1}y_i\frac{\partial
f }{\partial y_i}+h_{2i}x_i\frac{\partial f}{\partial x_i})$$ for
all elements $f\in A_n$. Note that each element of $K$ acts on
$A_n$ by a Poisson derivation.
\end{defn}

Let $A$ be a Poisson algebra and let an additive group $H$ act on $A$ by Poisson
derivations. A proper Poisson ideal $Q$ of $A$ is said to be $H$-prime
Poisson ideal if $Q$ is $H$-stable such that whenever $I,J$ are
$H$-stable Poisson ideals of $A$ with $IJ\subseteq Q$, either
$I\subseteq Q$ or $J\subseteq Q$. A Poisson algebra $A$ is said to
be $H$-simple if $0$ and $A$ are the only $H$-stable Poisson
ideals of $A$.

\begin{lem}
Let $A$ be a Poisson algebra and let $\alpha$ be a Poisson
derivation on $A$. Suppose that $H$ acts on $A[x^{\pm1};\alpha]_p$
so that $x$ is an $H$-eigenvector and $A$ is both $H$-stable and
$H$-simple, where $H$ acts on $A$ by restriction. If $H$ contains
a Poisson derivation $g$ such that $g|_A=\alpha$ and $g(x)=cx$
for some $0\neq c\in{\bf k}$ then $A[x^{\pm1};\alpha]_p$ is
$H$-simple.
\end{lem}
\begin{proof}
Let $I$ be a nonzero proper $H$-Poisson ideal of
$A[x^{\pm1};\alpha]_p$. Then choose $0\neq a\in I$, of shortest
length with respect to $x$, say $a=a_kx^k+\cdots+a_mx^m$ for some
$k\leq m$, where $a_i\in A$ for each $i$ and $a_k,a_m\neq0$. Since
$x$ is unit and $A\cap I=0$, we may assume that $k=0$ and
$a=a_0+\cdots+a_mx^m$, where $m>0$ and $a_0,a_m\neq0$. Set
$J=\{r\in A\ |\ r+r_1x+\cdots+r_mx^m\in I\mbox{ for some }
r_1,\cdots,r_m\in A\}$ and note that $J$ is a Poisson ideal of
$A$. Given any $h\in H$, let $\lambda_h$ be the $h$-eigenvalue of
$x$. Since $I$ is $H$-stable,
$h(r+r_1x+\cdots+r_mx^m)=h(r)+(h(r_1)+\lambda_hr_1)x+\cdots+(h(r_m)+m\lambda_hr_m)x^m\in
I$, and so $h(r)\in J$. Hence $J$ is an $H$-Poisson ideal of $A$, and
thus either $J=0$ or $J=A$; by our choice of $a$, $1\in J$. Thus
we may assume that $a=1+a_1x+\cdots, a_mx^m$. Since $I$ is
$H$-stable, $g(a)=(g(a_1)+ca_1)x+\cdots+(g(a_m)+mca_m)x^m\in I$, which
has the length less than $a$, hence $g(a)=0$ and
$g(a_i)+ica_i=\alpha(a_i)+ica_i=0$ for each $i=1,\cdots,m$. Now,
$\{a,x\}=\alpha(a_1)x^2+\cdots+\alpha(a_m)x^{m+1}$ is an element
of $I$ with the length less than $a$. Hence $\alpha(a_i)=0$ and thus
$a_i=0$ for all $i=1,\cdots,m$. It follows that $a=1\in I$, a contradiction.
As a result, $A[x^{\pm1};\alpha]_p$ is
$H$-simple.
\end{proof}

\begin{lem}
Let $B=A[y;\alpha]_p[x;\beta]_p$, where $A$ is a prime Poisson
algebra and both $\alpha$ and $\beta$ are Poisson derivations,
such that $\beta(A)\subseteq A$ and $\beta(y)=cy$ for some
$c\in{\bf k}$, and that $H$ is a group of Poisson derivations on
$B$ such that $A$ is $H$-stable and $y,x$ are $H$-eigenvectors. If
there exist $f,g\in H$ such that $f|_A=\alpha$ with $f(y)=ay$ and
$g|_{A[y;\alpha]_p}=\beta$ with $g(x)=bx$ for some $a,b\in{\bf
k}^\times$, and if $A$ is $H$-simple, then
\begin{itemize}
\item[(a)] $B[y^{-1}][x^{-1}],B/\langle y,x\rangle,(B/\langle
y\rangle)[x^{-1}]$, and $(B/\langle x\rangle)[y^{-1}]$ are
$H$-simple.
\item[(b)] $B$ has only four $H$-prime Poisson ideals $0,\langle
y\rangle,\langle x\rangle, \langle y,x\rangle$.
\end{itemize}
\end{lem}
\begin{proof}
(a) Note that
$$B[y^{-1}]=A[y^{\pm1};\alpha]_p[x;\beta]_p,\ \ \
B[y^{-1}][x^{-1}]=A[y^{\pm1};\alpha]_p[x^{\pm1};\beta]_p.$$ By
Lemma 3.2, $A[y^{\pm1};\alpha]_p$ is $H$-simple. Now apply Lemma
3.2 twice to obtain that
$B[y^{-1}][x^{-1}]=A[y^{\pm1};\alpha]_p[x^{\pm1};\beta]_p$ is
$H$-simple.

Since $B/\langle y,x\rangle\cong_H A$, it follows that $B/\langle
y,x\rangle$ is $H$-simple. Next, the Poisson algebra $(B/\langle
y\rangle)[x^{-1}]\cong_HA[x^{\pm1};\beta]_p$ is $H$-simple by
Lemma 3.2. Analogously, $(B/\langle
x\rangle)[y^{-1}]\cong_HA[y^{\pm1};\alpha]_p$ is $H$-simple.

(b) Clearly, $0,\langle y\rangle,\langle x\rangle, \langle
y,x\rangle$ are all $H$-prime Poisson ideals. Suppose that $P$ is
a nonzero $H$-prime Poisson ideal of $B$. The extended ideal
$P^e=PB[y^{-1}][x^{-1}]$ contains the multiplicative identy because
$B[y^{-1}][x^{-1}]$ is $H$-simple. Thus, $y^ix^j\in P$ for some
$i,j$ and thus $P$ contains $y$ or $x$ since $\langle y\rangle$
and $\langle x\rangle$ are both $H$-stable Poisson ideals of $B$.
If $x\in P$ then $P/\langle x\rangle$ is an $H$-prime Poisson
ideal of $B/\langle x\rangle$, and thus $P=\langle x\rangle$ or
$P=\langle x, y\rangle$ since $(B/\langle x\rangle)[y^{-1}]$ is
$H$-simple. Analogously, if $P$ contains $y$ then $P=\langle
y\rangle$ or $P=\langle x, y\rangle$. As a result, $B$ has only
four $H$-prime Poisson ideals $0,\langle y\rangle,\langle
x\rangle, \langle y,x\rangle$.
\end{proof}

\begin{lem}
Let $B=(A;\alpha,\beta,c,u)=A[y;\alpha]_p[x;\beta',\delta]_p$ be the Poisson algebra given in Proposition 1.4.
Assume, in addition, that
$A$ is a prime Poisson algebra, $\alpha(u)=du$, $\beta(u)=-du$ for some
$d\in{\bf k}$  with $c+d\neq0$ and $0\neq\delta(y)=u\in A$ is Poisson normal in $B$.
Set $z=(c+d)yx+\delta(y)$. Let $H$ be a group
of Poisson derivations on $B$ such that $A$ is $H$-stable and
$y,x$ and $z$ are $H$-eigenvectors. Suppose that there exist
$f,g\in H$ such that $f|_A=\alpha$ with $f(y)=ay$
for some $a\in{\bf k}^\times$ and $g|_{A[y;\alpha]_p}=\beta'$ with
$g(y^{-1}z)=by^{-1}z$ for some $b\in{\bf k^\times}$. If $A$ is $H$-simple,
then
\begin{itemize}
\item[(a)] $\delta(y)$ is invertible in $B$.
\item[(b)] no proper $H$-stable Poisson ideal of $B$ contains a
power of $y$.
\item[(c)] $B[y^{-1}][z^{-1}]$, $B[z^{-1}]$ and $B/\langle
z\rangle$ are $H$-simple.
\item[(d)] the only $H$-prime Poisson ideals of $B$ are $0$ and
$\langle z\rangle$.
\end{itemize}
\end{lem}
\begin{proof}
(a) Since $\delta(y)=\{y,x\}-cyx$ is $H$-eigenvector and Poisson
normal, $\langle\delta(y)\rangle$ is an $H$-stable Poisson ideal
of $B$. Thus $I=\langle\delta(y)\rangle\cap A$ is a nonzero
$H$-stable Poisson ideal of $A$, and hence $1\in I$ since $A$ is
$H$-simple. In particular, $1\in \langle\delta(y)\rangle$ and so
$\delta(y)B=\langle\delta(y)\rangle=B$. Consequently, $\delta(y)$
is invertible in $B$.

(b) Suppose that $P$ is a proper $H$-Poisson ideal of  $B$ such
that $y^j\in P$ for some $j>0$. Whenever $y^j\in P$ for some
$j>0$, we have that
$$jy^{j-1}\delta(y)=\delta(y^j)=\{y^j,x\}-\beta'(y^j)x=\{y^j,x\}-jcy^jx\in
P,$$ and hence $y^{j-1}\in P$ since $\delta(y)$ is invertible in
$B$ by (a). The repeated applications of the above argument
guarantee that $y\in P$. Therefore $\delta(y)=\{y,x\}-cyx\in P$,
and thus no proper $H$-Poisson ideal contains a power of $y$ since
$\delta(y)$ is invertible in $B$ by (a).

(c) By Lemma 1.7 (b), $B[y^{-1}]=A[y^{\pm1};\alpha]_p[y^{-1}z;\beta']_p$.
 Note that
$$B[y^{-1}][z^{-1}]=A[y^{\pm1};\alpha]_p[(y^{-1}z)^{\pm1};\beta']_p, \ \ \
g|_{A[y^{\pm1};\alpha]_p}=\beta'.$$
Applying Lemma 3.2 yields that both $A[y^{\pm1};\alpha]_p$ and
$A[y^{\pm1};\alpha]_p[(y^{-1}z)^{\pm1};\beta']_p$ are $H$-simple, so
$B[y^{-1}][z^{-1}]$ is $H$-simple.

Let $P$ be  an $H$-prime Poisson ideal of $B[z^{-1}]$. Then $P$ is
induced from an $H$-prime Poisson ideal $\check{ P}$ of $B$ disjoint
from $\{z^j\ |\ j=0,1,\cdots\}$. By (b), $\check P$ contains no
$y^j$. Suppose that $\check P$ contains some $y^iz^j$. Since $z$
and $y$ are Poisson normal and $H$-eigenvectors by Lemma 1.7 and the hypothesis, we have that
$y^i\in\check P$ or $z^j\in\check P$, a contradiction. Thus
$\check P$ is disjoint from the multiplicative set generated by
$y$ and $z$. Hence the extension $\check P^e$ to
$B[y^{-1}][z^{-1}]$ is an $H$-prime Poisson ideal. Since
$B[y^{-1}][z^{-1}]$ is $H$-simple, $\check P^e=0$, and so $\check
P=0$, so $P=0$. Thus $B[z^{-1}]$ contains no nonzero $H$-prime
Poisson ideals.

If $I$ is a proper $H$-Poisson ideal of $B[z^{-1}]$ then $I$ is
contained in a prime Poisson ideal $P$ of $B[z^{-1}]$ by Lemma 1.8.
Set $Q=(P:H)$ the largest $H$-stable Poisson ideal contained in $P$.
If $I$ and $J$ are $H$-stable Poisson ideals such that $IJ\subseteq Q$
then either $I\subseteq P$ or $J\subseteq P$, and thus either $I\subseteq Q$ or $J\subseteq Q$.
It follows that $Q$
is an $H$-prime Poisson ideal such that $I\subseteq
Q\subseteq P$. Since $B[z^{-1}]$ does not have a nonzero $H$-prime
Poisson ideal, we have that $I=Q=0$. Hence, $B[z^{-1}]$ is
$H$-simple.

Note that $\langle z\rangle$ is a Poisson ideal of $B$ since $z$
is Poisson normal by Lemma 1.7, and $zB[y^{-1}]$ is also a Poisson
ideal of $B[y^{-1}]$. Observe that
$$\aligned
(B/\langle z\rangle)[y^{-1}]&\cong_H B[y^{-1}]/(zB[y^{-1}])\\
&=A[y^{\pm1};\alpha]_p[y^{-1}z;\beta']_p/(zA[y^{\pm1};\alpha]_p[y^{-1}z;\beta']_p)\\
&\cong_H A[y^{\pm1};\alpha]_p .
\endaligned$$
 Thus $(B/\langle z\rangle)[y^{-1}]$ is
$H$-simple by Lemma 3.2. Denote by $\overline b$ the canonical homomorphic
image of $b\in B$ in $B/\langle z\rangle$. Since
$\overline{yx}=-(c+d)^{-1}\overline{\delta(y)}$ and $\delta(y)$ is
invertible in A by (a), $\overline y$ is invertible in $B/\langle
z\rangle$, and thus $B/\langle z\rangle=(B/\langle
z\rangle)[y^{-1}]$ is $H$-simple.

(d) Clearly $0$ is an $H$-prime Poisson ideal of $B$ since $B$ is
a prime Poisson algebra. Further, $\langle z\rangle$ is $H$-stable and prime
Poisson since $z$ is an $H$-eigenvector and Poisson normal in $B$.
Now, let $P$ be an $H$-prime Poisson ideal of $B$. If $P$ contains
no $z^i$ then $P$ extends to an $H$-prime Poisson ideal $\check P$
of $B[z^{-1}]$. Since $B[z^{-1}]$ is $H$-simple by (c), $\check P=0$, and
so $P=0$. Assume that $P$ contains some $z^i$. Then $z\in P$ since
$\langle z\rangle$ is an $H$-stable Poisson ideal and $P$ is an
$H$-prime Poisson ideal. Thus $0$ and $\langle z\rangle$ are the
only $H$-prime Poisson ideals of $B$ since $B/\langle z\rangle$
is $H$-simple by (c).
\end{proof}

\begin{defn}
Given an admissible set $T$ of $A_n$, let $N_T$ be the subset of
$\mathcal P_n$ defined by
\begin{itemize}
\item[(a)]
$y_1\in N_T$ if and only if $y_1\notin T$
\item[(b)] $x_1\in N_T$ if and only if $x_1\notin T$
\item[(c)] for $i>1$, $\Omega_i\in N_T$ if and only if $\Omega_{i-1}\notin T$ and $\Omega_i\notin T$
\item[(d)] for $i>1$, $y_i\in N_T$ if and only if $\Omega_{i-1}\in T$ and $y_i\notin T$
\item[(e)] for $i>1$, $x_i\in
N_T$ if and only if $\Omega_{i-1}\in T$ and $x_i\notin T$
\end{itemize}
\end{defn}

\begin{thm}
For an admissible set $T$, let $E_T$ be the multiplicative set
generated by $N_T$.
\begin{itemize}
\item[(a)] $E_T\cap \langle T\rangle=\phi$.
\item[(b)] $A_n^T=(A_n/\langle T\rangle)[E_T^{-1}]$ is $H$-simple.
\end{itemize}
\end{thm}
\begin{proof}
(a) It follows immediately from Lemma 2.6 (b).

(b)
We proceed by induction on $n$. Let  $n=1$ and we will apply Lemma
3.3(a). By Remark 2.2, $A_1=({\bf k},0,0, -q_1,0)={\bf
k}[y_1;0]_p[x_1,\beta_1]_p$, where $\beta_1(y_1)=-q_1y_1$, and
consider $f=(1,1), g=(-q_1,1)\in K$. Then $g$ acts as $\beta_1$ on
$A_1$ and $g(x)=x$.  There are four possible cases for $T$:
$$\phi, \{y_1, \Omega_1\}, \{x_1, \Omega_1\}, \{y_1, x_1, \Omega_1\}.$$
 Hence
$A_1^T$ is  one of the forms $A_1[y^{-1}][x^{-1}]$, $(A_1/\langle
y_1\rangle)[x_1^{-1}]$, $(A_1/\langle x_1\rangle)[y_1^{-1}]$,
$A_1/\langle y_1, x_1\rangle$. Applying Lemma 3.3(a), $A_1^T$ is
$H$-simple.

Suppose that $n>1$ and $A_{n-1}^S$ is $K$-simple for any
admissible set $S\subseteq \mathcal P_{n-1}$. Note that
$$\aligned
&A_n=A_{n-1}[y_n;\alpha_n]_p[x_n;\beta_n,\delta_n]_p=(A_{n-1};\alpha_n,\beta_n,-q_n,-\Omega_{n-1})\\
&\alpha_n(-\Omega_{n-1})=p_n(-\Omega_{n-1}), \ \ \ \beta_{n}(-\Omega_{n-1})=-p_n(-\Omega_{n-1})
\endaligned$$
by Remark 2.2 and Lemma 2.3.
Given an
admissible set $T$ of $A_n$, set $T'=T\cap\mathcal P_{n-1}$ and
let $I$ be the ideal of $A_{n-1}$ generated by $T'$. Then, since
$I$ is $\{\alpha_n,\beta_n,\delta_n\}$-stable, Lemma
1.3 gives the following $K$-equivalence:
$$A_n/IA_{n}\cong_K(A_{n-1}/I)[y_n;\overline\alpha_n]_p[x_n;\overline\beta_n,\overline\delta_n]_p,$$
where $\overline\delta_n=0$ if $\Omega_{n-1}\in T'$, and thus we have
$$(A_n/IA_{n})[E_{T'}^{-1}]\cong_K
(A_{n-1}/I)[E_{T'}^{-1}][y_n;\overline\alpha_n]_p[x_n;\overline\beta_n,\overline\delta_n]_p.$$

Set $A=(A_{n-1}/I)[E_{T'}^{-1}]$ and $S=T\backslash T'$. Then
$\langle T\rangle=IA_n+\langle S\rangle$ and
$$\aligned
A_n/\langle T\rangle&\cong_K (A_n/IA_n)/(\langle T\rangle/IA_n)\\
A_n/\langle T\rangle[E_{T'}^{-1}]&\cong_K
A[y_n;\overline\alpha_n]_p[x_n;\overline\beta_n,\overline\delta_n]_p/\langle
S\rangle.\endaligned$$
Let $E$ be the multiplicative set generated
by $N_T\backslash (N_{T'}\cap\mathcal P_{n-1})$. Then
$$\aligned
A_n/\langle T\rangle[E_{T}^{-1}]&=A_n/\langle T\rangle[E_{T'}^{-1}][E^{-1}]\\
&\cong_K
(A[y_n;\overline\alpha_n]_p[x_n;\overline\beta_n,\overline\delta_n]_p/\langle
S\rangle)[E^{-1}].\endaligned$$

In order to apply Lemma 3.3 and Lemma 3.4, we will define the
necessary elements of $K$. Set
$$\aligned
f&=(\gamma_{1n},p_n-\gamma_{1n},
\gamma_{2n},p_n-\gamma_{2n},\cdots,
\gamma_{n-1,n},p_n-\gamma_{n-1,n}, 1, p_n-1)\\
g&=(-q_1-\gamma_{1n}, q_1-p_n+\gamma_{1n}, -q_2-\gamma_{2n},
q_2-p_n+\gamma_{2n},\cdots, \\
&\qquad\qquad\qquad-q_{n-1}-\gamma_{n-1,n},
q_{n-1}-p_n+\gamma_{n-1,n},-q_n, q_n-p_n).
\endaligned$$
Then $f,g\in K$ and $f|_{A_{n-1}}=\alpha_n, f(y_n)=y_n,
f(x_n)=(p_n-1)x_n$ and $g|_{A_{n-1}[y_n;\alpha_n]_p}=\beta_n,
g(x_n)=(q_n-p_n)x_n$. Note that
$(-q_n+p_n)y_nx_n-\Omega_{n-1}=-\Omega_n$ and
$g(-y_n^{-1}\Omega_n)=(q_n-p_n)(-y_n^{-1}\Omega_n)$.
As defined, $1$ and $q_n-p_n$ are nonzero.

There are five possible cases for $S$:
$$\phi, \ \ \{\Omega_n\},\ \
\{y_n,\Omega_n\}, \ \ \{x_n,\Omega_n\},\ \  \{y_n,x_n,\Omega_n\}.$$
 If $S=\phi$ then
$\langle S\rangle=0$, and if $\Omega_{n-1}\in T'$, then $E$ is
generated by $y_n$ and $x_n$, so that
$$\aligned
A_n^T&\cong_K(A[y_n;\overline\alpha_n]_p[x_n;\overline\beta_n,\overline\delta_n]_p/\langle
S\rangle)[E^{-1}]\\
&=(A[y_n;\overline\alpha_n]_p[x_n;\overline\beta_n]_p)[y_n^{-1}][x_n^{-1}]\\
&=A[y_n^{\pm1};\overline\alpha_n]_p[x_n^{\pm1};\overline\beta_n]_p.
\endaligned$$
 since $\overline\delta_n=0$. Applying Lemma 3.3 yields that
$A_n^T$ is $K$-simple. If $\Omega_{n-1}\notin T'$ then $E$ is generated
by $\Omega_n$ and
$A_n^T\cong_K(A[y_n;\overline\alpha_n]_p[x_n;\overline\beta_n,\overline\delta_n]_p)[\Omega_n^{-1}]
$ is $K$-simple by Lemma 3.4.

If $S=\{\Omega_n\}$ then $\Omega_{n-1}\notin T'$ and $E=\{1\}$, and so
$$A_n^T\cong_K(A[y_n;\overline\alpha_n]_p[x_n;\overline\beta_n,\overline\delta_n]_p)/\langle
\Omega_n\rangle$$ is $K$-simple by Lemma 3.4.

If $S=\{y_n,\Omega_n\}$ then $\langle S\rangle=\langle y_n\rangle$ and
$E$ is generated by $x_n$. Further $\overline\delta_n=0$ since
$\Omega_{n-1}\in T'$ and
$$\aligned
A_n^T&\cong_K(A[y_n;\overline\alpha_n]_p[x_n;\overline\beta_n,\overline\delta_n]_p/\langle
S\rangle)[x_n^{-1}]\\
&=(A[y_n;\overline\alpha_n]_p[x_n;\overline\beta_n]_p/\langle
y_n\rangle)[x_n^{-1}]
\endaligned
$$
is $K$-simple by Lemma 3.3.

If $S=\{x_n,\Omega_n\}$ then $\langle S\rangle=\langle x_n\rangle$ and
$E$ is generated by $y_n$. Moreover $\overline\delta_n=0$ and
$$\aligned
A_n^T&\cong_K(A[y_n;\overline\alpha_n]_p[x_n;\overline\beta_n,\overline\delta_n]_p/\langle
S\rangle)[y_n^{-1}]\\
&=(A[y_n;\overline\alpha_n]_p[x_n;\overline\beta_n]_p/\langle
x_n\rangle)[y_n^{-1}]
\endaligned$$
is $K$-simple by Lemma 3.3.

Lastly, if $S=\{y_n,x_n,\Omega_n\}$ then $\Omega_{n-1}\in T'$ and $E=\{1\}$,
and so
$$A_n^T\cong_K(A[y_n;\overline\alpha_n]_p[x_n;\overline\beta_n]_p)/\langle
y_n,x_n\rangle$$ is $K$-simple by Lemma 3.3. Therefore we conclude
that $A_n^T$ is $K$-simple for every admissible set $T$.
\end{proof}

\begin{lem}
Let $P$ be a $K$-prime Poisson ideal of $A_n$. Then
$T=P\cap\mathcal P_n$ is an admissible set.
\end{lem}
\begin{proof}
For convenience, set $\Omega_0=0$. Suppose that $y_i\in T$,
$i=1,\cdots,n$. Then $\Omega_{i-1}=\{x_i,y_i\}-q_iy_ix_i\in P$ and
$\Omega_i=(q_i-p_i)y_ix_i+\Omega_{i-1}\in P$ by Lemma 2.3. It follows that if $y_i\in T$
then $\Omega_i,\Omega_{i-1}\in T$. Similarly, if $x_i\in T$, $i=1,\cdots,n$
then $\Omega_i,\Omega_{i-1}\in T$. Conversely, suppose that $\Omega_i,\Omega_{i-1}\in
T$, $i=1,\cdots,n$. Then $(q_i-p_i)y_ix_i= \Omega_i-\Omega_{i-1}\in P$.
Since $y_i$ and $x_i$ are both $K$-eigenvectors and Poisson normal
modulo $\langle \Omega_{i-1}\rangle$, we have that
$\langle y_i,\Omega_{i-1}\rangle$ and $\langle x_i,\Omega_{i-1}\rangle$
are $K$-stable Poisson ideals and
$\langle y_i,\Omega_{i-1}\rangle$$\langle x_i,\Omega_{i-1}\rangle\subseteq P,$
and hence we have $y_i\in P$ or $x_i\in P$. Therefore,
if $\Omega_i,\Omega_{i-1}\in T$, $i=1,\cdots,n$ then $y_i\in T$ or
$x_i\in T$. It follows that $T$ is an admissible set of $A_n$.
\end{proof}

\begin{thm}
Every $K$-prime Poisson ideal of $A_n$ is generated by an
admissible set.
\end{thm}
\begin{proof}
Let $P$ be a $K$-prime Poisson ideal of $A_n$ and let
$T=P\cap\mathcal P_n$. Then $T$ is an admissible set by Lemma 3.7
and $P/\langle T\rangle$ is a $K$-prime Poisson ideal of
$A_n/\langle T\rangle$. By definition, $N_T\cap T=\phi$ and so
$N_T\cap P=\phi$, and hence $\overline N_T\cap P/\langle
T\rangle=\phi$, where each element of $\overline N_T$ is Poisson
normal in $A_n/\langle T\rangle$. Recalling that $\overline E_T$
is the multiplicative set generated by $\overline N_T$, we have
that $\overline E_T\cap P/\langle T\rangle=\phi$. Hence
$(P/\langle T\rangle)[{{\overline E}_T}^{-1}]$ is a $K$-prime Poisson ideal
of $A_n^T$, and so $(P/\langle T\rangle)[{{\overline E}_T}^{-1}]=0$
since $A_n^T$ is $K$-simple by Theorem 3.6. Therefore, $P/\langle
T\rangle=0$, so $P=\langle T\rangle$.
\end{proof}

\section{A relation between Poisson spaces and Quantum spaces}

 We set $\Lambda=\Bbb Z^n$, the free abelian group
with finite rank $n$ and let $u : \Lambda\times
\Lambda:\longrightarrow {\bf k}$ be  an antisymmetric biadditive map, that is,
$$
u(\alpha,\beta)=-u(\beta,\alpha), \ \ \ u(\alpha,\beta+\beta')=u(\alpha,\beta)+u(\alpha,\beta'),$$
for $\alpha,\beta,\beta'\in \Lambda$.
Then the group algebra $k\Lambda$
becomes a Poisson algebra with Poisson bracket
$$\{x_\alpha,x_\beta\}=u(\alpha,\beta)x_{\alpha+\beta}.$$
This Poisson algebra is denoted by ${\bf k}_u\Lambda$. (See \cite[\S2]{Oh4}.)
Set $\Lambda^+=(\Bbb Z^+)^n$ and let ${\bf k}_u\Lambda^+$
be the subalgebra of ${\bf k}_u\Lambda$ generated by all $x_\alpha, \alpha\in\Lambda^+$.
Note that ${\bf k}_u\Lambda^+$ is isomorphic to the commutative polynomial ring
${\bf k}[x_1,\cdots,x_n]$. Hence ${\bf k}[x_1,\cdots,x_n]$ is a Poisson algebra with
Poisson bracket
$$\{x_i,x_j\}=u(\epsilon_i,\epsilon_j)x_ix_j$$
 for all $i,j$, where $\epsilon_i$ are the standard basis for $\Lambda$.

Let ${\bf r}=(r_{ij})$ be a skew-symmetric $n\times n$-matrix with
entries in ${\bf k}$. Then there exists an antisymmetric biadditive map
$u : \Lambda\times
\Lambda:\longrightarrow {\bf k}$
defined by $u(\epsilon_i,\epsilon_j)=r_{ij}$. Hence
${\bf k}[x_1,\cdots,x_n]$ is a Poisson algebra with
Poisson bracket
$$\{x_i,x_j\}=r_{ij}x_ix_j$$
for all $i,j$. This Poisson algebra is denoted by $Q({\bf r})$.

\begin{defn}
Let ${\bf r}_n$ be the following skew-symmetric $2n\times
2n$-matrix determined by the defining coefficients in
$A_n=A_{n,\Gamma}^{P,Q}$ given in Theorem 2.1:\footnotesize
$$ \bordermatrix{
        &Y_1   &X_1    &Y_2          &X_2               &\cdots   &Y_n         &X_n   \cr
    Y_1 &0     &-q_1   &\gamma_{12}  &-(q_1+\gamma_{12}) &\cdots   &\gamma_{1n}  &-(q_1+\gamma_{1n}) \cr
    X_1 &q_1   &0      &p_2-\gamma_{12} &q_1-p_2+\gamma_{12}&\cdots &p_n-\gamma_{1n} &q_1-p_n+\gamma_{1n} \cr
    Y_2 &-\gamma_{12}&-(p_2-\gamma_{12}) &0 &-q_2       &\cdots  &\gamma_{2n}   &-(q_2+\gamma_{2n}) \cr
    X_2 &q_1+\gamma_{12}&-(q_1-p_2+\gamma_{12})&q_2 &0    &\cdots  &p_n-\gamma_{2n} &q_2-p_n+\gamma_{2n} \cr
    \vdots &\vdots &\vdots &\vdots &\vdots &\ddots &\vdots       &\vdots \cr
    Y_n &-\gamma_{1n}&-(p_n-\gamma_{1n}) &-\gamma_{2n} &-(p_n-\gamma_{2n}) &\cdots &0         &-q_n \cr
    X_n &q_1+\gamma_{1n} &-(q_1-p_n+\gamma_{1n})&q_2+\gamma_{2n}&-(q_2-p_n+\gamma_{2n}) &\cdots &q_n  &0
}$$\normalsize Hence there exists the Poisson algebra $Q({\bf
r}_n)={\bf k}[Y_1,X_1,\cdots,Y_n,X_n]$, which is called  the
Poisson algebra attached to $A_n$.
\end{defn}

\begin{defn} \cite[Definition 1.1]{Hor1}
Let $P,Q\in{({\bf k})^\times}^n$ such that $P=(p_1,\cdots,p_n)$
and $Q=(q_1,\cdots,q_n)$ where $p_iq_i^{-1}$ is not a root of
unity for each $i=1,\cdots,n$. Further, let
$\Gamma=(\gamma_{ij})\in M_n({\bf k}^\times)$ be a multiplicative
skew-symmetric matrix, that is, $\gamma_{ji}=\gamma_{ij}^{-1}$ and
$\gamma_{ii}=1$ for all $i,j$. Then $K_{n,\Gamma}^{P,Q}$ is
the ${\bf k}$-algebra generated by $y_1,x_1,\cdots, y_n,x_n$
satisfying the following relations:
\[\begin{array}{ll}
y_iy_j=\gamma_{ij}y_jy_i &(\mbox{all } i,j)\\
x_iy_j=p_j\gamma_{ij}^{-1}y_jx_i  &(i<j)\\
y_ix_j=q_i^{-1}\gamma_{ij}^{-1}x_jy_i &(i<j)\\
x_ix_j=q_ip_j^{-1}\gamma_{ij}x_jx_i  &(i<j)\\
x_iy_i=q_iy_ix_i+\sum_{k=1}^{i-1}(q_k-p_k)y_kx_k &(\mbox{all } i)
 \end{array}\]
 We simply write $K_n$ for $K_{n,\Gamma}^{P,Q}$ unless any confusion arises.

 As in \cite[Lemma 2.1 and Definition 2.2]{Hor1}, set
$$\Omega_{i}=\sum_{k=1}^{i}(q_k-p_k)y_kx_k\in K_n$$ for each $i=1,\cdots
,n$ and a subset $T\subseteq\{y_1,x_1,\Omega_1,\cdots,
y_n,x_n,\Omega_n\}=\mathcal{P}_n\subseteq K_n$ is said to be  an admissible set
of $K_n$ if $T$ satisfies the conditions (a) and (b) of
Definition 2.4.
\end{defn}

Let $\bold s=(s_{ij})$ be a multiplicative
skew-symmetric $n\times n$ matrix over $k$. Then the multiparameter
quantized coordinate ring $R(\bold s)$ of affine
$n$-space is the algebra generated by $x_1,\cdots,x_n$ subject to
the relations
$$x_ix_j=s_{ij}x_jx_i.$$

\begin{defn} \cite[(5)]{GoToKa}
Let ${\bf s}_n$ be the following multiplicative skew-symmetric
$2n\times 2n$-matrix determined by the defining coefficients in
$K_n=K_{n,\Gamma}^{P,Q}$ given in Definition 4.2:
$$ \bordermatrix{
        &Y_1   &X_1    &Y_2          &X_2               &\cdots   &Y_n         &X_n   \cr
    Y_1 &1     &q^{-1}_1   &\gamma_{12}  &q_1^{-1}\gamma_{12}^{-1} &\cdots   &\gamma_{1n}  &q_1^{-1}\gamma_{1n}^{-1} \cr
    X_1 &q_1   &1  &p_2\gamma_{12}^{-1} &q_1p_2^{-1}\gamma_{12}&\cdots &p_n\gamma_{1n}^{-1} &q_1p_n^{-1}\gamma_{1n} \cr
    Y_2 &\gamma_{12}^{-1}&p_2^{-1}\gamma_{12} &1 &q_2^{-1} &\cdots &\gamma_{2n} &q_2^{-1}\gamma_{2n}^{-1} \cr
    X_2 &q_1\gamma_{12}&q_1^{-1}p_2\gamma_{12}^{-1}&q_2 &1 &\cdots &p_n\gamma_{2n}^{-1} &q_2p_n^{-1}\gamma_{2n} \cr
    \vdots &\vdots &\vdots &\vdots &\vdots &\ddots &\vdots       &\vdots \cr
    Y_n &\gamma_{1n}^{-1}&p_n^{-1}\gamma_{1n} &\gamma_{2n}^{-1} &p_n^{-1}\gamma_{2n} &\cdots &1 &q_n^{-1} \cr
    X_n &q_1\gamma_{1n} &q_1^{-1}p_n\gamma_{1n}^{-1}&q_2\gamma_{2n}&q_2^{-1}p_n\gamma_{2n}^{-1} &\cdots &q_n
    &1
}$$

Note that each entry of ${\bf r}_n$ given in Definition 4.1 is
the additive form of the corresponding entry of ${\bf s}_n$. Since
${\bf s}_n$ is a multiplicative skew-symmetric matrix, there
exists the multiparameter quantized coordinate ring $R({\bf s}_n)$ of affine
$2n$-space,  which is called the quantized algebra attached to $K_n$.
\end{defn}

Let $I$ be an ideal of a Poisson algebra $A$. Denote by $\mathcal{V}_p(I)$ the set of all
prime Poisson ideals of $A$ containing $I$. That is,
$$\mathcal{V}_p(I)=\mathcal{V}(I)\cap\pspec A.$$
Since $J=\cap\mathcal{V}_p(I)$ is a semiprime Poisson ideal and $\mathcal{V}_p(I)=\mathcal{V}_p(J)$,
the closed sets of $\pspec A$ are exactly the sets $\mathcal{V}_p(I)$ for semiprime Poisson ideals $I$ of $A$.

\begin{lem}
Let $A$ be a finitely generated Poisson algebra.
Then the irreducible closed sets of $\pspec A$ are exactly the sets
$\mathcal{V}_p(P)$ for prime Poisson ideals $P$.
\end{lem}
\begin{proof}
Let $P$ be a prime Poisson ideal of $A$ and let
$\mathcal{V}_p(P)=\mathcal{V}_p(I)\cup\mathcal{V}_p(J)$ for some ideals $I,J$.
Then we have
$$P=\cap \mathcal{V}_p(P)=\cap[\mathcal{V}_p(I)\cup\mathcal{V}_p(J)]=
\cap\mathcal{V}_p(I\cap J)\supseteq I\cap J\supseteq IJ,$$
and thus $P$ contains $I$ or $J$, say $I\subseteq P$. Hence we have that
$\mathcal{V}_p(P)=\mathcal{V}_p(I)$. It follows that $\mathcal{V}_p(P)$ is irreducible.

Conversely, let $\mathcal{V}_p(P)$ be an irreducible set for some semiprime Poisson ideal $P$. Since $A$
is noetherian and $P$ is semiprime, $P=Q_1\cap\cdots\cap Q_n$ an intersection of finitely many prime ideals minimal
over $P$ by \cite[Theorem 2.4]{GoWa}. Since each $(Q_i:\mathcal{H}(A))$ is a prime Poisson ideal
containing $P$ by \cite[Lemma 1.3]{Oh4} and $Q_i$ is a prime minimal over $P$, we have
$Q_i=(Q_i:\mathcal{H}(A))$  a prime Poisson ideal. Since
$\mathcal{V}_p(P)=\mathcal{V}_p(Q_1\cap\cdots\cap Q_n)=\cup_i\mathcal{V}_p(Q_i)$ and $\mathcal{V}_p(P)$ is irreducible,
we have that  $\mathcal{V}_p(P)=\mathcal{V}_p(Q_i)$ for some $i$. It completes the proof.
\end{proof}

\begin{defn}
Given an admissible set $T$ of $A_n$ (or $K_n$), set
$$\aligned
&\eta(T)=\{Y_i\ |\ y_i\in T\}\cup\{X_i\ |\ x_i\in T\}\cup\{X_i\ |\
\Omega_i\in T, y_i\notin T, x_i\notin T\}\\
&\mathcal Y_T \mbox{ the multiplicative set of $A_n$ (or $K_n$)
generated by $\{y_i|\
y_i\notin T\}$}\\
&\mathcal U_T \mbox{ the multiplicative set of $Q({\bf r}_n)$ (or
$R({\bf s}_n)$)
generated by  $\{Y_i|\ y_i\notin T\}$}.
\endaligned$$

Note that if $\Omega_i\in T$ then $Y_iX_i\in\langle\eta(T)\rangle$  since
$Y_i\in \eta(T)$ or $X_i\in \eta(T)$.
\end{defn}

\begin{lem}
For admissible sets $T$ and $T'$, $T\neq T'$ $\Leftrightarrow$
$\eta(T)\neq\eta(T')$.
\end{lem}
\begin{proof}
($\Rightarrow$) Let $T\neq T'$. If there exists $\Omega_i$ such
that $\Omega_i\in T$ and $\Omega_i\notin T'$, or $\Omega_i\notin
T$ and $\Omega_i\in T'$, say $\Omega_i\in T$ and $\Omega_i\notin
T'$, then $Y_i\in \eta(T)$ or $X_i\in\eta(T)$, but $\eta(T')$
contains neither $Y_i$ nor $X_i$, and so $\eta(T)\neq\eta(T')$.
Hence we suppose that
$T\cap\{\Omega_1,\cdots,\Omega_n\}=T'\cap\{\Omega_1,\cdots,\Omega_n\}$.
If $T\cap\{y_1,\cdots,y_n\}\neq T'\cap\{y_1,\cdots,y_n\}$ then
$\eta(T)\neq\eta(T')$ clearly, and if $T\cap\{y_1,\cdots,y_n\}=
T'\cap\{y_1,\cdots,y_n\}$ then there exists  $x_i$ such that one
of $T$ and $T'$ contains $x_i$ and the other does not contain
$x_i$, say $x_i\in T$ and $x_i\notin T'$.  If $i=1$ then
$X_1\in\eta(T)$ and $X_1\notin\eta(T')$, and so
$\eta(T)\neq\eta(T')$. Assume that $i>1$. Then $\Omega_{i-1}\in T$
and so $\Omega_{i-1}\in T'$, and hence $y_i\in T'$ since
$\Omega_i\in T'$ and $x_i\notin T'$.  It follows that
$X_i\in\eta(T)$ and $X_i\notin\eta(T')$, and thus
$\eta(T)\neq\eta(T')$.

($\Leftarrow$) If $T=T'$ then $\eta(T)=\eta(T')$ clearly.
Therefore, if $\eta(T)\neq\eta(T')$ then we have that $T\neq T'$.
\end{proof}

\begin{lem} Let $T$ be an admissible set of $A_n$.

(a) The map $\Psi'_T:(A_n/\langle T\rangle)[{\overline{\mathcal
Y}_T}^{-1}]\longrightarrow (Q({\bf r}_n)/\langle
\eta(T)\rangle)[{\overline{\mathcal U}_T}^{-1}]$ defined by
$$\begin{array}{ll}
\Psi'_T(\overline y_i)=\overline Y_i &(\mbox{all $i$})\\
\Psi'_T(\overline x_1)=\overline X_1 &(i=1) \\
\Psi'_T(\overline x_i)=\overline
X_i-(\hat{q}_i-\hat{p}_i){\overline Y_i}^{-1}\overline
Y_{i-1}\overline X_{i-1} &(i\geq2,y_i\notin T, x_i\notin T,\Omega_{i-1}\notin T)\\
\Psi'_T(\overline x_i)=\overline X_i &(i\geq2, x_i\notin T, \Omega_{i-1}\in T)\\
\Psi'_T(\overline x_i)=-(\hat{q}_i-\hat{p}_i){\overline
Y_i}^{-1}\overline Y_{i-1}\overline X_{i-1} &(i\geq2, y_i\notin T,
x_i\notin T, \Omega_i\in T)\end{array}$$
 is a Poisson isomorphism, where $(\hat{q}_i-\hat{p}_i)=(q_i-p_i)^{-1}(q_{i-1}-p_{i-1})$
 for each $i=2,\cdots,n$.
 Moreover $\Psi'_T(\overline{\mathcal
Y}_T)=\overline{\mathcal U}_T$, $\Psi'_T(\overline \Omega_i)=(q_i-p_i)\overline{Y}_i\overline{X}_i$
for all $i$ and if $f\in A_n$ and $T'$ is an admissible set
such that $T\subseteq T'$ then $\Psi'_T(f)$
is congruent to $\Psi'_{T'}(f)$ modulo the ideal of $(Q({\bf r}_n)/\langle
\eta(T')\rangle)[{\overline{\mathcal U}_{T'}}^{-1}]$ generated by $\eta(T')$.

 (b)
Set
$$\aligned
\spec_TA_n&=\{P\in\spec A_n\ |\ P\cap\mathcal{P}_n=T\}\\
\spec_{\eta(T)}Q({\bf r}_n)&=\{P\in\spec Q({\bf r}_n)\ |\ P\cap\{Y_1,X_1,\cdots, Y_n,X_n\}=\eta(T)\}\\
\pspec_{\eta(T)}Q({\bf r}_n)&=\{P\in\pspec Q({\bf r}_n)\ |\ P\cap\{Y_1,X_1,\cdots, Y_n,X_n\}=\eta(T)\}.
\endaligned$$
 For every $P\in\text{spec}_T A_n$, there exists a unique element
 $\Psi_T(P)$ of $\text{spec}_{\eta(T)}Q({\bf r}_n)$ such that
  $\Psi'_T((P/\langle T\rangle)^e)=(\Psi_T(P)/\langle \eta(T)\rangle)^e$, and conversely, for
each $Q\in\spec_{\eta(T)}Q({\bf r}_n)$, there exists a unique $P\in\spec_TA_n$ such that
$\Psi'_T((P/\langle T\rangle)^e)=(Q/\langle \eta(T)\rangle)^e$.
  Moreover, the map
  $$\Psi_T:\text{pspec}_TA_n\longrightarrow \text{pspec}_{\eta(T)}Q({\bf r}_n),\ \  P\mapsto\Psi_T(P)$$
  is a homeomorphism, if $P$ and $Q$ are prime ideals of $A_n$ such that
$P\subseteq Q$, $P\cap\mathcal{P}_n=T, Q\cap\mathcal{P}_n=T'$ then $T\subseteq T'$
and $\Psi_T(P)\subseteq\Psi_{T'}(Q)$,
  and $P\in\text{pspec}_TA_n$ is symplectic if and only if $\Psi_T(P)$ is symplectic.

  (c) The map $\Psi:\text{pspec}A_n=\biguplus\text{pspec}_T
  A_n\longrightarrow \biguplus \text{pspec}_{\eta(T)}Q({\bf r}_n)$
  defined by $\Psi(P)=\Psi_T(P)$ for $P\in \text{pspec}_T
  A_n$ is a homeomorphism such that its restriction
  $\Psi|_{\text{symp} A_n}:\text{symp} A_n\longrightarrow \Psi(\text{symp} A_n)$ is also a homeomorphism.
\end{lem}
\begin{proof}
(a) Define an algebra homomorphism
$$\psi_T:A_n\longrightarrow
(Q({\bf r}_n)/\langle \eta(T)\rangle)[{\overline{\mathcal
U}_T}^{-1}]$$ by
\begin{equation}\begin{array}{ll}
\psi_T(y_i)=\overline Y_i &(\text{all } i)\\
\psi_T(x_i)=\overline X_i-(\hat{q}_i-\hat{p}_i){\overline
Y_i}^{-1}\overline
Y_{i-1}\overline X_{i-1} &(i\geq2, y_i\notin T)\\
\psi_T(x_i)=\overline X_i &(i\geq2, y_i\in T)\\
\psi_T(x_1)=\overline X_1 &(i=1)\end{array}
\end{equation} where
$(\hat{q}_i-\hat{p}_i)=(q_i-p_i)^{-1}(q_{i-1}-p_{i-1})$ for each
$i=2,\cdots,n$. It is checked by a straight calculation that $\psi_T$ is a Poisson homomorphism.

We will prove that $\psi_T(\Omega_i)=(q_i-p_i)\overline Y_i\overline
X_i$ for all $i=1,\cdots,n$ using the induction on $i$.
If $i=1$ then $\psi(\Omega_1)=(q_1-p_1)\overline Y_1\overline
X_1$ clearly. Suppose now that $i>1$ and that
$\psi_T(\Omega_j)=(q_j-p_j)\overline Y_j\overline
X_j$ for each $j<i$. If $y_i\notin T$ then, by the induction hypothesis,
$$\aligned
\psi_T(\Omega_i)&=\psi_T((q_i-p_i)y_ix_i+\Omega_{i-1})\\
&=(q_i-p_i)\overline Y_i(\overline X_i-(\hat{q}_i-\hat{p}_i){\overline
Y_i}^{-1}\overline
Y_{i-1}\overline X_{i-1})+(q_{i-1}-p_{i-1})\overline Y_{i-1}\overline
X_{i-1}\\
&=(q_i-p_i)\overline Y_i\overline X_i.
\endaligned$$
If $y_i\in T$ then $\Omega_{i-1}\in T$ and thus
$Y_{i-1}X_{i-1}\in\langle\eta(T)\rangle$ and
$$\aligned
\psi_T(\Omega_i)&=\psi_T((q_i-p_i)y_ix_i+\Omega_{i-1})\\
&=(q_i-p_i)\overline Y_i\overline X_i+(q_{i-1}-p_{i-1})\overline Y_{i-1}\overline
X_{i-1}\\
&=(q_i-p_i)\overline Y_i\overline X_i.
\endaligned$$
Hence $\psi_T(\Omega_i)=(q_i-p_i)\overline Y_i\overline
X_i$ for all $i$, as required.

If $\Omega_i\in T$ then
$Y_iX_i\in\langle\eta(T)\rangle$, and so $\psi_T(\Omega_i)=0$ for all $\Omega_i\in T$.
Moreover, $\psi_T(x_i)=0$ for all $x_i\in T$ and $\psi_T(y_i)=0$
for all $y_i\in T$ clearly. It follows that $\langle
T\rangle\subseteq \ker \psi_T$ and thus there exists a Poisson
epimorphism $\overline{ \psi}_T:(A_n/\langle
T\rangle)[{\overline{\mathcal Y}_T}^{-1}]\longrightarrow (Q({\bf
r}_n)/\langle \eta(T)\rangle)[{\overline{\mathcal U}_T}^{-1}]$
such that $\overline{\psi}_T(\overline f)=\psi_T(f)$ for all $f\in
A_n$ since all elements of $\psi_T(\mathcal Y_T)$ are invertible.

Let us show that $\overline{\psi}_T=\Psi'_T$. Clearly
$\overline{\psi}_T(\overline y_i)=\Psi'_T(\overline y_i)$ for all
$i$ and $\overline{\psi}_T(\overline x_1)=\Psi'_T(\overline x_1)$.
Moreover it is clear that $\overline{\psi}_T(\overline
x_i)=\Psi'_T(\overline x_i)$ for all $x_i$ such that $y_i\notin T$
and $\Omega_{i-1}\notin T$, where $i\geq2$. If $i\geq2$ and
$y_i\in T$ then $\Omega_{i-1}\in T$ and so
$\overline{\psi}_T(\overline x_i)=\overline X_i=\Psi'_T(\overline
x_i)$. If $i\geq2$, $y_i\notin T, x_i\notin T$ and $\Omega_i\in T$
then
$\psi_T((q_i-p_i)y_ix_i)=\psi_T(\Omega_i-\Omega_{i-1})=-(q_{i-1}-p_{i-1})\overline
Y_{i-1}\overline X_{i-1}$ and thus $\overline{\psi}_T(\overline
x_i)=-(\hat{q}_i-\hat{p}_i){\overline{Y}_i}^{-1}\overline
Y_{i-1}\overline X_{i-1}=\Psi'_T(\overline x_i)$. Thus
$\overline{\psi}_T(\overline x_i)=\Psi'_T(\overline x_i)$ for all
$x_i$, as required.

Let $T'$ be an admissible set such that $T\subseteq T'$. Then
$\Psi'_T(y_i)$ and $\Psi'_T(x_i)$ are congruent to $\Psi'_{T'}(y_i)$ and $\Psi'_{T'}(x_i)$, respectively,
modulo the ideal of $(Q({\bf r}_n)/\langle
\eta(T')\rangle)[{\overline{\mathcal U}_{T'}}^{-1}]$ generated by $\eta(T')$ since
$\Psi'_{T'}(\Omega_i)=(q_i-p_i)\overline Y_i\overline X_i$, and thus
$\Psi'_T(f)$ is congruent to $\Psi'_{T'}(f)$ modulo the ideal  generated by $\eta(T')$ for all $f\in A_n$.

Note that the number of elements in $\eta(T)$ is equal to
$\text{length}(T)$ and thus the Gelfand-Kirillov dimension of
$(Q({\bf r}_n)/\langle \eta(T)\rangle)[{\overline{\mathcal
U}_T}^{-1}]$ is equal to $2n-\text{length}(T)$ by \cite[Example
3.6 and Proposition 4.2]{KrLe}. By Lemma 2.10 and \cite[Proposition
4.2]{KrLe},  $(A_n/\langle T\rangle)[{\overline{\mathcal
Y}_T}^{-1}]$ has the  Gelfand-Kirillov dimension
$2n-\text{length}(T)$. Moreover $\ker \overline\psi_T$ is a prime ideal since
$(Q({\bf
r}_n)/\langle \eta(T)\rangle)[{\overline{\mathcal U}_T}^{-1}]$ is
prime, and thus $\overline\psi_T=\Psi'_T$  is a Poisson
isomorphism by \cite[Proposition 3.16]{KrLe}.

(b) For every $P\in\text{spec}_T A_n$, the extended ideal
$(P/\langle T\rangle)^e$ is a prime ideal of $(A_n/\langle
T\rangle)[{\overline{\mathcal Y}_T}^{-1}]$ since $P\cap \mathcal
Y_T=\phi$, and thus, by (a), there exists a prime ideal $\Psi_T(P)$ of $Q({\bf
r}_n)$ such that
\begin{equation}\aligned
&\langle\eta(T)\rangle\subseteq\Psi_T(P)\\
 &\Psi_T(P)\cap\mathcal{U}_T=\phi\\
&\Psi'_T((P/\langle T\rangle)^e)=(\Psi_T(P)/\langle \eta(T)\rangle)^e.
\endaligned
\end{equation}

If $Y_i\notin\eta(T)$ then $Y_i\notin\Psi_T(P)$ by the second formula of (4.2), and thus
if $Y_i\in\Psi_T(P)$ then $Y_i\in\eta(T)$. Let
$X_i\in\Psi_T(P)-\eta(T)$. Then  we have $\Psi'_T(\overline{\Omega}_i)=(q_i-p_i)\overline{Y}_i\overline{X}_i\in
(\Psi_T(P)/\langle \eta(T)\rangle)^e$ by (a) and thus $\Omega_i\in P$ by the third formula of (4.2).
It follows that $\Omega_i\in T$ and $x_i\notin T$ since $X_i\notin\eta(T)$. If $i=1$ then
$\Psi'_T(\overline x_1)=\overline X_1$ and if $i>1$ then $y_i\in T$ by the definition of $\eta(T)$ since
$\Omega_i\in T,x_i\notin T, X_i\notin\eta(T)$, and so $\Omega_{i-1}\in T$. Thus
$\Psi'_T(\overline x_i)=\overline X_i\in(\Psi_T(P)/\langle\eta(T)\rangle)^e$ for any $i$ by the definition of
$\Psi'_T$. Hence
$x_i\in P$ by the third formula of (4.2) and thus $x_i\in T$ and $X_i\in\eta(T)$, a contradiction. Therefore we have that
$\Psi_T(P)\cap\{Y_1,X_1,\cdots,Y_n,X_n\}=\eta(T)$ by the first formula of (4.2).
Conversely if $Q\in\spec_{\eta(T)}Q({\bf r}_n)$ then there exists a unique element $P\in\spec_TA_n$
such that $\Psi'_T((P/\langle T\rangle)^e)=(Q/\langle \eta(T)\rangle)^e$ since $\Psi_T'$
is an algebra isomorphism.

Since $\Psi'_T$ is a Poisson isomorphism, we have that
$P\in\pspec_TA_n$ if and only if $\Psi_T(P)\in\pspec_{\eta(T)}Q({\bf r}_n)$ and the map
$$\Psi_T:\text{pspec}_TA_n\longrightarrow \text{pspec}_{\eta(T)}Q({\bf r}_n),\ \  P\mapsto\Psi_T(P)$$
is a homeomorphism.
Let $P$ and $Q$ be prime ideals of $A_n$ such that
$P\subseteq Q$, $P\cap\mathcal{P}_n=T$ and $Q\cap\mathcal{P}_n=T'$.
Then $T$ and $T'$ are admissible sets such that $T\subseteq T'$.
Since, for each $f\in A_n$,
$\Psi'_T(f)$ is congruent to $\Psi'_{T'}(f)$ modulo the ideal
of $Q({\bf r}_n)[\mathcal{U}_{T'}^{-1}]$ generated by $\eta(T')$ by (a), we have that
$\Psi_T(P)\subseteq\Psi_{T'}(Q)$.

The fact that $P\in\text{pspec}_TA_n$ is symplectic if and only if $\Psi_T(P)$ is symplectic
will be proved after the proof of Theorem 4.10.

(c) By (b), it is clear that $\Psi$ is bijective. For a prime Poisson ideal $P$ of $A_n$ with
$P\cap\mathcal{P}_n=T$,
let us show that
\begin{equation}
\Psi(\mathcal{V}_p(P))=\mathcal{V}_p(\Psi_T(P))\cap \text{Im}\Psi.
\end{equation}
If (4.3) is true then $\Psi^{-1}$ is continuous by Lemma 4.4.
Let $Q\in\mathcal{V}_p(P)$ and $Q\cap\mathcal{P}_n=T'$. Then $P\subseteq Q$ and thus
$\Psi_T(P)\subseteq \Psi_{T'}(Q)=\Psi(Q)$ by (b). Conversely,
let $Q$ be a prime Poisson ideal of $A_n$ such that $Q\cap\mathcal{P}_n=T'$
and $\Psi_{T'}(Q)\in\mathcal{V}_p(\Psi_T(P))$.
Then $\langle\eta(T)\rangle\subseteq\Psi_T(P)\subseteq\Psi_{T'}(Q)$ and thus
$T\subset Q$. It follows that $T\subseteq T'$ and $\Psi_T(P)\subseteq\Psi_{T'}(Q)$.
Therefore $P\subseteq Q$ since
$\Psi'_T(f)$ is
congruent to $\Psi'_{T'}(f)$ modulo the ideal of $Q({\bf r}_n)[\mathcal{U}_{T'}^{-1}]$ generated by $\eta(T')$
 for all $f\in A_n$, and hence
$\Psi(\mathcal{V}_p(P))=\mathcal{V}_p(\Psi_T(P))\cap \text{Im}\Psi$, as claimed.

For any ideal $I$ of $Q({\bf r}_n)$, we have that
$\mathcal{V}_p(I)\cap\text{Im}\Psi=\bigcup_T[\mathcal{V}_p(I)\cap\text{Im}\Psi_T]$ and
each set $\mathcal{V}_p(I)\cap\text{Im}\Psi_T$
is a union of the forms $\mathcal{V}_p(Q)\cap\text{Im}\Psi_T$ for some prime Poisson ideals $Q$ of
$Q({\bf r}_n)$ with  $Q\cap\{Y_1,X_1,\cdots,Y_n,X_n\}=\eta(T)$ by (b) and Lemma 4.4. Hence if we show that
$\Psi^{-1}(\mathcal{V}_p(Q))=\mathcal{V}_p(\Psi_T^{-1}(Q))$
for a prime Poisson ideal $Q$ of $Q({\bf r}_n)$ with $Q\cap\{Y_1,X_1,\cdots,Y_n,X_n\}=\eta(T)$
then $\Psi$ is continuous. By (4.3), we have that
$$\mathcal{V}_p(Q)\cap\text{Im}\Psi=\Psi(\mathcal{V}_p(\Psi_T^{-1}(Q)))$$
for a prime Poisson ideal $Q$ of $Q({\bf r}_n)$ with $Q\cap\{Y_1,X_1,\cdots,Y_n,X_n\}=\eta(T)$,
and thus we have that $\Psi^{-1}(\mathcal{V}_p(Q))=\mathcal{V}_p(\Psi_T^{-1}(Q))$.
It completes the proof.
\end{proof}

For a Poisson algebra $A$, the set
$$Z_p(A)=\{a\in A\ |\ \{a,A\}=0\}$$ is said to be a Poisson center of $A$.

\begin{defn} \cite[Theorem 2.4]{Oh4}
A Poisson ${\bf k}$-algebra $A$ is said to satisfy the Poisson Dixmier-Moeglin equivalence if the following
conditions are equivalent: For a prime Poisson ideal $P$ of $A$,
\begin{itemize}
\item[(i)] $P$ is symplectic (i.e., there exists a maximal ideal $M$ of $A$ such that $P$ is the largest
Poisson ideal contained in $M$).
\item[(ii)] $P$ is rational (i.e., the Poisson center of the quotient field of $A/P$ is algebraic over ${\bf k}$).
\item[(iii)] $P$ is locally closed (i.e., the intersection of all prime Poisson ideals properly containing $P$
is strictly larger than $P$).
\end{itemize}
\end{defn}

\begin{prop} Let ${\bf r}$ be a skew-symmetric $n\times n$-matrix with entries in ${\bf k}$.
The Poisson algebra $Q({\bf r})$ satisfies the Poisson Dixmier-Moeglin equivalence.
More precisely,
$$\aligned
\symp Q({\bf r})&=\{\text{locally closed prime Poisson ideals}\}\\
&=\{\text{rational prime Poisson ideals}\}\\
&=\bigcup_W\{\text{maximal elements of }\pspec_WQ({\bf r})\},
\endaligned$$
where $W\subseteq\{x_1,x_2,\cdots,x_n\}$ and
$$\pspec_WQ({\bf r})=\{P\in\pspec Q({\bf r})\ |\ P\cap\{x_1,x_2,\cdots,x_n\}=W\}.$$
\end{prop}
\begin{proof}
Let $P$ be a prime Poisson ideal of $Q({\bf r})$.
If $P$ is locally closed then $P$ is symplectic by \cite[Proposition 1.7]{Oh4}
since $Q({\bf r})$ is a Jacobson ring, and
if $P$ is symplectic then $P$ is rational by \cite[Proposition 1.10]{Oh4}.

Let $P$ be a rational prime ideal and set $W=P\cap\{x_1,x_2,\cdots,x_n\}$. Then
$Q({\bf r})/\langle W\rangle\cong {\bf k}[z_1,\cdots,z_k]$, where
$\{z_1,\cdots,z_k\}=\{x_1,x_2,\cdots, x_n\}-W$, and thus
$\pspec_WQ({\bf r})$ is homeomorphic to $\pspec {\bf k}[z_1,\cdots,z_k]$ and
$P$ corresponds to a rational
prime Poisson ideal $Q$ of ${\bf k}[z_1^{\pm1},\cdots,z_k^{\pm1}]$.
Since $Q$ is a maximal prime Poisson ideal of ${\bf k}[z_1^{\pm1},\cdots,z_k^{\pm1}]$
by \cite[Corollary 2.3  and Theorem 2.4]{Oh4}, $P$ is a maximal element of $\pspec_WQ({\bf r})$.

If $P$ is a maximal element of $\pspec_WQ({\bf r})$ and $Q$ is a prime Poisson ideal of $Q({\bf r})$
which  contains $P$ properly then $Q$ has an element $x_i$ which is not in $P$. Hence the intersection
of all prime Poisson ideals properly containing $P$ is strictly larger than $P$. It follows that
$P$ is locally closed.
\end{proof}

\begin{thm}
The Poisson algebra $A_n$ satisfies the Poisson Dixmier-Moeglin equivalence.
More precisely,
$$\aligned
\symp A_n&=\{\text{locally closed prime Poisson ideals}\}\\
&=\{\text{rational prime Poisson ideals}\}\\
&=\bigcup_T\{\text{maximal elements of }\pspec_TA_n\}
\endaligned$$
\end{thm}
\begin{proof}
Let $P$ be a prime Poisson ideal of $A_n$. Note that $A_n$ is a Jacobson ring since $A_n$ is
a finitely generated commutative ring.
If $P$ is locally closed then $P$ is symplectic by \cite[Proposition 1.7]{Oh4}, and
if $P$ is symplectic then $P$ is rational by \cite[Proposition 1.10]{Oh4}.

Let $P$ be rational and set $T=P\cap \mathcal P_n$. Then $T$ is an admissible set by Proposition 2.8 and
$\Psi_T(P)$ is a rational prime Poisson ideal of $Q({\bf r}_n)$ such that
$\Psi_T(P)\cap\{Y_1,X_1,\cdots,Y_n,X_n\}=\eta(T)$ by Lemma 4.7 (b). Hence
 $\Psi_T(P)$ is a maximal element of
the set $\pspec_{\eta(T)} Q({\bf r}_n)$ by Proposition 4.9,
and thus $P$ is a maximal element of $\pspec_TA_n$ by Lemma 4.7 (b).

Let $P$ be a maximal element of $\pspec_TA_n$ and let $Q$ be a prime Poisson ideal of $A_n$ which
properly contains $P$. Then $Q\cap\mathcal P_n$ contains properly the admissible set $T$. Hence
 the intersection of all prime Poisson ideals
$Q$ of $A_n$ properly containing $P$
is strictly larger than $P$ since admissible sets are finite. It follows
that $P$ is locally closed. It completes the proof.
\end{proof}

{\it The proof of Lemma 4.7 (b)}:
For any $P\in\text{pspec}_TA_n$, $P$ is symplectic if and only if $P$ is
a maximal element of $\text{pspec}_TA_n$ by Theorem 4.10, and $\Psi_T(P)$
is a symplectic ideal of $Q({\bf r}_n)$ if and only if $\Psi_T(P)$
is a maximal element of $\text{pspec}_{\eta(T)}Q({\bf r}_n)$ by Proposition 4.9.
Hence $P$ is symplectic
if and only if $\Psi_T(P)$ is symplectic since $\Psi_T$ preserves inclusion by Lemma 4.7 (b).
 $\square$

\medskip
The map $\Upsilon'_T$ of  the following Lemma 4.11 (a) is modified from the map
given in \cite[Proposition 2.7]{GoToKa}.

\begin{lem} Let $T$ be an admissible set of $K_n=K_{n,\Gamma}^{P,Q}$.

(a) The map $\Upsilon'_T:(K_{n,\Gamma}^{P,Q}/\langle
T\rangle)[{\overline{\mathcal Y}_T}^{-1}]\longrightarrow (R({\bf
s}_n)/\langle \eta(T)\rangle)[{\overline{\mathcal U}_T}^{-1}]$
defined by
$$\begin{array}{ll}
\Upsilon'_T(\overline y_i)=\overline Y_i &(\mbox{all $i$})\\
\Upsilon'_T(\overline x_1)=\overline X_1 &(i=1) \\
\Upsilon'_T(\overline x_i)=\overline
X_i-(\hat{q}_i-\hat{p}_i){\overline Y_i}^{-1}\overline
Y_{i-1}\overline X_{i-1} &(i\geq2,y_i\notin T, x_i\notin T,\Omega_{i-1}\notin T)\\
\Upsilon'_T(\overline x_i)=\overline X_i &(i\geq2, x_i\notin T, \Omega_{i-1}\in T)\\
\Upsilon'_T(\overline x_i)=-(\hat{q}_i-\hat{p}_i){\overline
Y_i}^{-1}\overline Y_{i-1}\overline X_{i-1} &(i\geq2,
y_i\notin T, x_i\notin T,\Omega_i\in T)\end{array}$$
 is an algebra isomorphism, where $(\hat{q}_i-\hat{p}_i)=(q_i-p_i)^{-1}(q_{i-1}-p_{i-1})$ for $i\geq2$.
 Moreover $\Upsilon'_T(\overline{\mathcal
Y}_T)=\overline{\mathcal U}_T$, $\Upsilon'_T(\overline \Omega_i)=(q_i-p_i)\overline{Y}_i\overline{X}_i$
for all $i$
and $\Upsilon'_T(\overline x_i)$ is congruent to $\overline X_i$ modulo
$\langle\overline{Y}_{i-1}\overline{X}_{i-1}\rangle$ for all $i=1,\cdots,n$, where
$\langle\overline{Y}_{0}\overline{X}_{0}\rangle=0$.

 (b)
 Let $\text{spec}_{\eta(T)}R({\bf s}_n)$ be the set of all prime
 ideals $Q$ of $R({\bf s}_n)$ such that
 $Q\cap\{Y_1,X_1,\cdots,Y_n,X_n\}=\eta(T)$.
 For every $P\in\text{spec}_T K_{n,\Gamma}^{P,Q}$, there exists a unique prime ideal
 $\Upsilon_T(P)\in\text{spec}_{\eta(T)}R({\bf s}_n)$ such that
  $\Upsilon'_T(P/\langle T\rangle)^e=(\Upsilon_T(P)/\langle \eta(T)\rangle)^e$.
  Moreover, the map $\Upsilon_T:\text{spec}_TK_{n,\Gamma}^{P,Q}\longrightarrow
  \text{spec}_{\eta(T)}R({\bf s}_n)$
  is a homeomorphism, if $P$ and $Q$ are prime ideals of  $K_{n,\Gamma}^{P,Q}$
  such that $P\subseteq Q$, $P\cap\mathcal{P}_n=T$, $Q\cap\mathcal{P}_n=T'$
  then $T\subseteq T'$ and $\Upsilon_T(P)\subseteq\Upsilon_{T'}(Q)$, and
   $P\in\text{spec}_T K_{n,\Gamma}^{P,Q}$
  is primitive if and only if $\Upsilon_T(P)$ is primitive.

  (c) The map $\Upsilon:\text{spec}K_{n,\Gamma}^{P,Q}=\biguplus_T\text{spec}_T
  K_{n,\Gamma}^{P,Q}\longrightarrow \biguplus \text{spec}_{\eta(T)}R({\bf s}_n)$
  defined by $\Upsilon(P)=\Upsilon_T(P)$ for $P\in \text{spec}_T
  K_{n,\Gamma}^{P,Q}$ is a homeomorphism such that its restriction
  $\Upsilon|_{\text{prim} K_{n,\Gamma}^{P,Q}}:\text{prim}K_{n,\Gamma}^{P,Q}\longrightarrow
  \Upsilon(\text{prim} K_{n,\Gamma}^{P,Q})$
  is also a homeomorphism.
\end{lem}
\begin{proof}
(a) It follows by mimicking the proofs of \cite[Proposition 2.7]{GoToKa} and Lemma 4.7(a).

(b) It follows immediately from mimicking the proof of Lemma 4.7(b).
Note that a prime ideal $P\in\spec_T K_n$ is primitive if and only if $P$ is a
maximal element of $\spec_T K_n$ by \cite[Theorem 3.8 and Theorem 4.16]{Hor1}.

(c) Note that the irreducible closed sets of $\spec K_n$ are exactly the sets
$\mathcal{V}(P)$ for prime ideals $P$ of $K_n$ and that $K_n$ is noetherian
by \cite[Theorem 4.12]{Hor1}.
 The result follows immediately from an argument mimicking the proof of Lemma 4.7 (c).
\end{proof}

\begin{rem}Let $P,Q,\Gamma$ be the ones given in Definition 4.2 and let $G$
be  the multiplicative subgroup of ${\bf k}^\times$ generated by
all $p_i, q_j,\gamma_{ij}$. Then, by \cite[Lemma 2.4]{OhPaSh1}, there exists a group homomorphism
$\phi$ from $G$ into ${\bf k}$ such that its restriction to the
torsion free subgroup of $G$ is injective. Since each
$p_iq_i^{-1}$ is not a root of unity, we have
$\phi(p_i)\neq\phi(q_i)$  for $i=1,2,\cdots, n$. Set
$$\aligned
\phi(P)&=(\phi(p_1),\cdots,\phi(p_n))\\
\phi(Q)&=(\phi(q_1),\cdots,\phi(q_n))\\
\phi(\Gamma)&=(\phi(\gamma_{ij})).
\endaligned$$

Note that $\phi(\Gamma)$ is a skew-symmetric $n\times n$-matrix with
entries in ${\bf k}$ since $\Gamma$ is a multiplicative skew-symmetric
$n\times n$-matrix.
\end{rem}

\begin{lem}
 Let $\pi_1:X\longrightarrow Y$ and  $\pi_2:Y\longrightarrow Z$
be topological quotient maps. Then the composition
$\pi_2\circ\pi_1:X\longrightarrow Z$ is  also a topological
quotient map.
\end{lem}
\begin{proof}
Clearly, the composition $\pi_2\circ\pi_1$ is surjective. For a subset
$C$ of $Z$, $C$ is closed in $Z$ if and only if $\pi^{-1}_2(C)$ is closed in $Y$
if and only if $\pi^{-1}_1(\pi^{-1}_2(C))=(\pi_2\circ\pi_1)^{-1}(C)$ is closed in $X$
since both $\pi_1$ and $\pi_2$ are topological quotient maps. Hence
$\pi_2\circ\pi_1$ is a topological quotient map.
\end{proof}

\begin{thm} Let $P,Q,\Gamma,\phi$ be the ones given in Definition
4.2 and Remark 4.12, and assume that the subgroup of ${\bf k}^\times$ generated by all $p_i,q_j,\gamma_{ij}$
does not contain $-1$.
Then there exists a topological quotient map
$$\pi:\text{pspec}A_{n,\phi(\Gamma)}^{\phi(P),\phi(Q)}
\longrightarrow\text{spec}K_{n,\Gamma}^{P,Q}$$ such that its restriction
$$\pi|_{\text{symp}A_{n,\phi(\Gamma)}^{\phi(P),\phi(Q)}}:\text{symp}A_{n,\phi(\Gamma)}^{\phi(P),\phi(Q)}
\longrightarrow \prim K_{n,\Gamma}^{P,Q}$$
is also a topological quotient map. Moreover, if $\phi$ is a monomorphism then both $\pi$ and
$\pi|_{\text{symp}A_{n,\phi(\Gamma)}^{\phi(P),\phi(Q)}}$ are homeomorphisms.
\end{thm}
\begin{proof}
By  Theorem 2.1 and Remark 4.12, the Poisson algebra
$A_{n,\phi(\Gamma)}^{\phi(P),\phi(Q)}$ is constructed. Set
$A=A_{n,\phi(\Gamma)}^{\phi(P),\phi(Q)}$. Let ${\bf
s}_n=(s_{ij})$ be the $2n\times 2n$-matrix attached to
$K_{n,\Gamma}^{P,Q}$ which is given in Definition 4.3. Note that $\phi({\bf
s}_n)=(\phi(s_{ij}))$ is the attached matrix to the Poisson
algebra $A$ by Definition 4.1 and Definition 4.3, and that $Q(\phi({\bf s}_n))$ is the Poisson
algebra attached to $A$.
Note that $Q(\phi({\bf s}_n))$ is Poisson isomorphic to ${\bf k}_u\Lambda^+$
for a suitable antisymmetric biadditive map, where $\Lambda^+=(\Bbb Z^+)^{2n}$,
that $R({\bf s}_n)$ is the multiparameter quantized coordinate ring of affine $2n$-space
corresponding to the antisymmetric $2n\times 2n$-matrix ${\bf s}_n$,
and that $\biguplus_T\text{pspec}_{\eta(T)}Q(\phi({\bf s}_n))\subseteq\pspec Q(\phi({\bf s}_n))$
and  $\biguplus_T\text{spec}_{\eta(T)}R({\bf s}_n)\subseteq\spec R({\bf s}_n)$.
Let
$$\aligned
&\pi_2:\pspec Q(\phi({\bf s}_n))\longrightarrow \spec R({\bf s}_n)\\
&\pi_2|_{\symp Q(\phi({\bf s}_n))}:\symp Q(\phi({\bf s}_n))\longrightarrow \prim R({\bf s}_n)
\endaligned$$
be the topological quotient maps given in
\cite[Corollary 3.6]{OhPaSh1}. Then the restriction of $\pi_2$
to the surjection
$$\aligned
\pi_1:\biguplus_T\text{pspec}_{\eta(T)}Q(\phi({\bf s}_n))&\longrightarrow
\pi_2(\biguplus_T\text{pspec}_{\eta(T)}Q(\phi({\bf s}_n)))=\biguplus_T\text{spec}_{\eta(T)}R({\bf s}_n)\\
P&\mapsto \pi_2(P)
\endaligned
$$
is a topological quotient map
by  \cite[Proof of Proposition 3.4]{OhPaSh1} and
 the restriction of $\pi_1$ to $\biguplus_T\text{symp}_{\eta(T)}Q(\phi({\bf s}_n))$
$$\pi_1|_{\biguplus_T\text{symp}_{\eta(T)}Q(\phi({\bf
s}_n))}:\biguplus_T\text{symp}_{\eta(T)}Q(\phi({\bf
s}_n))\longrightarrow\biguplus_T\text{prim}_{\eta(T)}R({\bf s}_n)
$$
is also a topological quotient map by \cite[Proof of Proposition 3.4]{OhPaSh1}
since the map $\pi_1|_{\biguplus_T\text{symp}_{\eta(T)}Q(\phi({\bf
s}_n))}$ is the restriction of $\pi_2|_{\symp Q(\phi({\bf s}_n))}$ to the set
${\biguplus_T\text{symp}_{\eta(T)}Q(\phi({\bf
s}_n))}$. Moreover,
 if $\phi$ is a monomorphism then both $\pi_1$ and its restriction $\pi_1|_{\biguplus_T\text{symp}_{\eta(T)}Q(\phi({\bf
s}_n))}$ are
homeomorphisms. Let $\Psi$ and $\Upsilon$ be the homeomorphisms given in Lemma 4.7 and Lemma 4.11,
respectively. Then the
composition $\pi=\Upsilon^{-1}\circ\pi_1\circ\Psi$
$$\aligned
\text{pspec}(A)=&\biguplus_T\text{pspec}_T(A){\overset\Psi\longrightarrow}
\biguplus_T\text{pspec}_{\eta(T)}Q(\phi({\bf r}_n))\\
&{\overset{\pi_1}\longrightarrow}\biguplus_T \text{spec}_{\eta(T)}R({\bf
r}_n){\overset{\Upsilon^{-1}}\longrightarrow}
\biguplus_T\text{spec}_TK_{n,\Gamma}^{P,Q}=\text{spec}
  K_{n,\Gamma}^{P,Q}\endaligned$$
is a topological quotient map such that its restriction to
$\text{symp}A$ is also a topological quotient map by  Lemma 4.13,
 and both $\pi$ and
$\pi|_{\text{symp}A}$ are homeomorphisms if the group homomorphism
$\phi$ is injective.
\end{proof}

\begin{lem}
Let $A$ be a finitely generated Poisson algebra. Then the map
$$\pi:\spec A\longrightarrow \pspec A, \ \ P\mapsto (P:\mathcal{H}(A))$$
is a topological quotient map.
\end{lem}
\begin{proof}
Since $\mathcal{H}(A)$ is the set of all Hamiltonians of $A$, if $P$ is a prime ideal then
$(P:\mathcal{H}(A))$ is a prime Poisson ideal of $A$ by \cite[Lemma 1.3]{Oh4}. Hence $\pspec A$ satisfies
the conditions (a) and (b) of \cite[1.1]{GoLet2} and the map $\pi$ is a topological quotient map
by \cite[1.7]{GoLet2}.
\end{proof}

\begin{lem}
Let $P$ be a symplectic ideal of $A_n$. Then
$$P=\cap\{M\in\max A_n\ |\ (M:\mathcal{H}(A_n))=P\}.$$
\end{lem}
\begin{proof}
By Proposition 2.8, the set $P\cap\mathcal{P}_n=T$ is an admissible set.
Let $M$ be a maximal ideal of $A_n$  such that $M\cap\mathcal{P}_n=T$ and $P\subseteq M$. Then
$(M:\mathcal{H}(A_n))=P$ since $(M:\mathcal{H}(A_n))$ is a symplectic ideal containing $P$ and $P$
is a maximal element of $\pspec_TA_n$ by Theorem 4.10.
Hence it is enough to show that there exists a set $U$ of maximal ideals $M$ of $A_n$ such that
$M\cap\mathcal{P}_n=T$, $P\subseteq M$ and $\cap U=P$.

Let $M$ be a prime ideal of $A_n$ such that $M\cap\mathcal{P}_n=T$ and
 $\Psi_T(M)$ is a
maximal ideal of $Q({\bf r}_n)$.
Then the extension $(M/\langle T\rangle)^e$  is a maximal ideal of
$(A_n/\langle T\rangle)[{\overline{\mathcal
Y}_T}^{-1}]$ and $\Psi_T(M)\cap\{Y_1,X_1,\cdots,Y_n,X_n\}=\eta(T)$ by Lemma 4.7 (b). If there exists a maximal ideal $N$ of $A_n$ properly containing $M$ then
$\Psi_T(M)\subset \Psi_{T'}(N)$ by Lemma 4.7 (b), where $T'=N\cap\mathcal{P}_n$. It is a
contradiction to the maximality of $\Psi_T(M)$. Hence if $M$ is a prime ideal  of $A_n$
such that $M\cap\mathcal{P}_n=T$ and
 $\Psi_T(M)$ is a maximal ideal of $Q({\bf r}_n)$
then $M$ is a maximal ideal of $A_n$.

The Poisson algebra $Q({\bf r}_n)/\langle\eta(T)\rangle$ is Poisson isomorphic to
$Q({\bf r}'_n)={\bf k}[z_1,\cdots,z_k]$,
where ${\bf r}'_n$ is the submatrix of ${\bf r}_n$ deleting rows and culumns corresponding to $Y_i$'s and
$X_i$'s in $\eta(T)$, and $k=2n-|\eta(T)|$. Set $B={\bf k}[z_1^{\pm1},\cdots,z_k^{\pm1}]$.
Note that if $N$ is a maximal ideal of $B$ then there exists a maximal ideal $M$ of
${\bf k}[z_1,\cdots,z_k]$ such that $M\cap\{z_1,\cdots,z_k\}=\phi$ and $MB=N$ and that
the symplectic ideal $\Psi_T(P)$
corresponds naturally to a maximal prime Poisson ideal $Q$ of $B$ by Proposition 4.9
or \cite[Corollary 2.3 and Theorem 2.4]{Oh4}.
Since $B={\bf k}[z_1^{\pm1},\cdots,z_k^{\pm1}]$ is a Jacobson ring, $Q$
is an intersection of maximal ideals of $B$.
Hence the symplectic ideal $\Psi_T(P)$ is an intersection of maximal ideals $N$ of $Q({\bf r}_n)$ such that
$\Psi_T(P)\subseteq N$ and $N\cap\{Y_1,X_1,\cdots, Y_n, X_n\}=\eta(T)$.
It follows by the above paragraph and Lemma 4.7 (b) that $P$ is
an intersection of maximal ideals $M$ of $A_n$ such that $P\subseteq M$ and $M\cap\mathcal{P}_n=T$.
It completes the proof.
\end{proof}

\begin{lem}
The map
$$\pi|_{\max A_n}:\max A_n\longrightarrow\symp A_n, \ M \mapsto(M:\mathcal{H}(A_n))$$
is a topological quotient map.
\end{lem}
\begin{proof}
The map
 $$\pi:\spec A_n\longrightarrow\pspec A_n, P\mapsto (P:\mathcal{H}(A_n))$$
is a topological map by Lemma 4.15 and its restriction
$$\pi|_{\max A_n}:\max A_n\longrightarrow\symp A_n, \ M \mapsto(M:\mathcal{H}(A_n))$$
is surjective by the definition of symplectic ideals.
  A subset $U$ of $\max A_n$ is said to be relatively $\pi$-stable
provided that $U=\max A_n\cap V$ for some subset $V\subseteq\spec A_n$ with  $V=\pi^{-1}(\pi(V))$.
If a subset $U$ of $\max A_n$ is relatively $\pi$-stable then $\cap U$ is an intersection of symplectic
ideals of $A_n$ by Lemma 4.16, and thus
the restriction $\pi|_{\max A_n}:\max A_n\longrightarrow \symp A_n$ is also a topological quotient map
by \cite[Lemma 1.3]{Oh4}.
\end{proof}

\begin{cor}
Let $P,Q,\Gamma$ be the ones given in Definition
4.2 and assume that the subgroup of ${\bf k}^\times$ generated by all $p_i,q_j,\gamma_{ij}$
does not contain $-1$.
Then there exists a topological quotient map
$$\pi:\spec {\bf k}[y_1,x_1,\cdots, y_n,x_n]\longrightarrow\spec  K_{n,\Gamma}^{P,Q}$$
 such that its restriction
$$\pi|_{\max {\bf k}[y_1,x_1,\cdots, y_n,x_n]}:\max {\bf k}[y_1,x_1,\cdots, y_n,x_n]
\longrightarrow\prim  K_{n,\Gamma}^{P,Q}$$ is also a topological
quotient map.
\end{cor}
\begin{proof}
It follows immediately from  Lemma 4.13, Theorem 4.14, Lemma 4.15 and Lemma 4.17.
\end{proof}

\bibliographystyle{amsplain}

\providecommand{\bysame}{\leavevmode\hbox to3em{\hrulefill}\thinspace}
\providecommand{\MR}{\relax\ifhmode\unskip\space\fi MR }
\providecommand{\MRhref}[2]{%
  \href{http://www.ams.org/mathscinet-getitem?mr=#1}{#2}
}
\providecommand{\href}[2]{#2}

\end{document}